\documentclass[reqno,oneside]{amsart}

\usepackage[utf8]{inputenc}

\usepackage[margin=3cm]{geometry}
\usepackage[dvipsnames]{xcolor}
\usepackage{latexsym}
\usepackage{mathrsfs}
\usepackage{graphicx}
\usepackage{amsmath,relsize,amssymb}
\usepackage{comment,enumerate,soul}
\usepackage{cite}
\usepackage{amsfonts}

\ifx\pdfoutput\undefined
       \usepackage[pdfstartview=FitB,pdfpagemode=UseNone,hypertex,colorlinks=true,linkcolor=red]{hyperref}
\else
       \usepackage[pdfstartview=FitB,pdfpagemode=UseNone,pdftex,hypertexnames=false,colorlinks=true,linkcolor=red]{hyperref}
\fi

\DeclareMathOperator*{\argmin}{arg\,min}

\usepackage[list=true]{subcaption}
\usepackage[sc,osf]{mathpazo}
\usepackage[foot]{amsaddr}

\usepackage{titlesec}

\titleformat{\section}[runin]
{\color{RoyalBlue}\bfseries}{\color{RoyalBlue}\thesection .}{1ex}{}[.]

\titleformat{\subsection}[runin]
{\color{RoyalBlue}\bfseries}{\color{RoyalBlue}\thesubsection .}{1ex}{}[.]

\titleformat{\subsubsection}[runin]
{\color{RoyalBlue}\bfseries}{\color{RoyalBlue}\thesubsubsection .}{1ex}{}[.]

\colorlet{ROYALBLUE}{RoyalBlue}

\usepackage{amsthm}

\newtheoremstyle{mystyle}%
{3pt}
{3pt}
{\itshape}
{}
{\bfseries\color{OliveGreen}}
{.}
{.5em}
{}

\theoremstyle{mystyle}
\newtheorem{remark}{Remark}

\newtheorem{theorem}{Theorem}[section]

\begin{document}

\title[\color{RoyalBlue}Ocean-depth measurement using SWW models]{\color{RoyalBlue}Ocean-depth measurement using shallow-water wave models}  
\author{\color{RoyalBlue}V Vasan}
\address[VV \& M]{International Centre for Theoretical Sciences, Tata Institute of Fundamental Research, India}
\email[VV]{vishal.vasan@icts.res.in}
\author{Manisha}


\author{\color{RoyalBlue}D Auroux}
\address[DA]{Universit\'e C\^ote d'Azur, CNRS, Laboratoire J.A. Dieudonn\'e, France }

\begin{abstract}
	In this paper, we consider a problem inspired by the real-world need to identify the topographical features of ocean basins. Specifically we consider the problem of estimating the bottom impermeable boundary to an inviscid, incompressible, irrotational fluid from measurements of the free-surface deviation alone, within the context of dispersive shallow-water wave models. The need to consider the shallow-water regime arises from the ill-posed nature of the problem and is motivated by prior work. We design an algorithm using which, both fluid velocities and the bottom-boundary profile, may be accurately recovered assuming an \textit{a priori} relatively inaccurate guess for the bottom boundary. We achieve this by considering two separate inverse problems: one to deduce the bottom-boundary from velocities and the surface deviation, and another to recover the velocities from the surface deviation and an approximate bottom-boundary. The former is a classic inverse problem that requires the inversion of an ill-conditioned matrix while the latter employs the observer framework. Combining the two inverse problems leads to our reconstruction algorithm. We emphasise the role played by model selection and its impact on algorithm design and the accuracy of the reconstruction.
	
	\bigskip

\end{abstract}

\maketitle




\section{Introduction}

The bottom of the ocean is as varied and diverse as the land above, replete with valleys, cliffs, seamounts and plateaus that are formed by a variety of geological processes \cite{brown2013ocean}. The precise shape of the bottom boundary determines ocean circulation and mixing which impacts Earth's climate \cite{kunze2004role} and the bio-diversity of the seas \cite{koslow1997seamounts}. Near-shore bathymetry is vital for the management of coastlines \cite{gorman1998monitoring} and to predict tsunami inundation \cite{yoon2002propagation,mori2011survey}. Determining the shape of the oceanic bottom boundary is not only of utmost importance, but also one of the most challenging oceanography problems both theoretically and practically. A common method to determine the bottom topography is via underwater acoustics using echo sounders \cite{MakrakisTaroudakis,CollinsKuperman}. However this method is slow, expensive and dangerous while leaving much of the shallow continental margin areas under-surveyed. Hence more modern approaches using satellite-gravity models that infer topography from deviations in the gravitational field have been promoted \cite{becker2009global}. These methods work on the principle that the additional mass due to a seamount increases the strength of the local gravitational field causing a bulge in the surface height of the water directly above \cite{smith1994bathymetric}. Despite the success of satellite-gravity based bathymetry \cite{smith2004conventional}, recent work suggests shipboard data is more reliable to detect sharp relief features with characteristic length-scales less than $25$km \cite{Watts2020}.

An alternative approach employs fluid dynamical principles to model the motion of an incompressible fluid, bounded below by a solid impermeable surface. With the rise in availability and quality of satellite imagery, fluid dynamical methods for ocean-depth measurement offer a way to survey large parts of the coastlines and oceans. The simplest of the fluid-mechanical methods consider variations in the dispersion relation of shoaling waves \cite{PiotrowskiDugan} and nonlinear corrections to these formulae \cite{Grilli}. All fluid dynamical methods are ultimately based on an analysis of the kinematic boundary condition at the free surface
\[\eta_t = w - \bar u_s \cdot \nabla\eta,\]
where $w$ and $\bar u_s$ are the vertical and horizontal velocities evaluated at the free surface $\eta$. The free surface $\eta$ and the surface velocities $(w,\bar u_s)$ are functionally dependent on the bottom boundary via the full equations of motion and associated boundary conditions. 

In the context of fully non-hydrostatic irrotational flow the kinematic boundary condition may be rewritten in terms of a surface potential $q$
\begin{align}
\eta_t = G(\eta,\zeta)q,\label{eqn:DNO}
\end{align}
where $G$ is the Dirichlet-Neumann operator (defined in Section \ref{sec:shallow-ww}) which, loosely speaking, maps the tangential fluid-velocity at the free surface to the normal fluid-velocity. This operator depends on the shape of the bottom-boundary profile $\zeta$. In \cite{nicholls2009detection} the authors considered the problem of estimating $G$ (and hence inferring $\zeta$) from frequencies and profiles of standing waves. Despite the ill-posed nature of the problem they were successful in determining the bottom-boundary from this information. On the other hand, Fontelos et al. \cite{fontelos2017bottom} assumed knowledge of $\eta,\eta_t,q$ at one instance of time and proved that there was a unique $\zeta$ that satisfied \eqref{eqn:DNO}. Moreover they proved that the solution to a minimisation problem allowed one to recover this bottom profile. Here too the authors noted the ill-posed nature of the problem which implies the reconstruction procedure is highly sensitive to noise in the input data. 

The defining property of the bottom-boundary is that the fluid-velocity normal to that surface vanishes. By employing a harmonic continuation argument,   Vasan \& Deconinck \cite{vasan2013inverse} wrote a nonlinear nonlocal equation which vanished at the true bottom-boundary. However this equation depended on the surface velocities. Although satellite measurements can infer both sea-surface deviations and sea-surface velocities, often the determination of the velocities requires an estimate of the bottom boundary \cite{chapron2005direct,johannessen2008direct}. As a result, determining the surface velocities as part of the bottom-boundary reconstruction, or at least assuming some characterisation of these velocities, seems inevitable. 

Equation \eqref{eqn:DNO} is only half the story. For a complete description of the fluid motion one must account for the momentum balance which, in the context of inviscid, irrotational, incompressible fluid flow, determines the time evolution of the surface potential $q$. Vasan \& Deconinck \cite{vasan2013inverse} used this equation for $q$ and the fact that the bottom-boundary was stationary in time to deduce a second nonlocal nonlinear equation, hence obtaining two equations for two unknowns: the bottom-boundary and the surface velocities. They were able to show reconstruction was possible based purely on surface data, \textit{i.e.} measurements of the surface deviation alone. Additionally they did not require \textit{a priori} knowledge of the mean depth. The authors noted the ill-posed nature of the problem and argued the reconstruction was more reliable in the shallow-water regime.

The Saint-Venant equations are a specific model for shallow-water waves widely used in modelling inundation of coastal regions. Upon averaging the incompressibility condition in the vertical and taking into account the kinematic boundary condition, one obtains
\begin{align*}
	\eta_t + \nabla\cdot ( (\eta  + h +  \zeta)\:\bar u) = 0,
\end{align*}
where $(\eta + h +  \zeta)$ represents the height of the fluid column and $\bar u$ is the depth-averaged horizontal velocity. Through a judicious choice of scaling, some time-averaging of the data and a knowledge of the mean-depth, the authors of \cite{kar2018ocean} considered the above equation as a hyperbolic PDE for the bottom-boundary profile with known coefficients. They assumed surface deviation and surface velocities were given by satellites and employed a small Rossby number parameter regime to determine the ageostrophic component of the velocity field. However the necessity of a non-zero Coriolis parameter limits the utility of this method for fully non-hydrostatic coastal regions. Although the fluid dynamical models considered thus far all involved long length-scales and inviscid flow, domain identification is not limited to such scenarios. Indeed one can pose the problem in the context of creeping flows \cite{heining2017flow}. 

Our approach closely follows the work of \cite{vasan2013inverse,fontelos2017bottom} but we view the problem, from the outset, in the shallow-water regime. Specifically we derive the shallow-water analogue of the equations in \cite{vasan2013inverse,fontelos2017bottom}. The idea to consider the shallow-water re-formulation was inspired by a conversation between Harvey Segur and one of the authors, several years ago. While VV attempted to describe the ill-posed nature of the problem, Harvey remarked the problem would `perhaps be easier, if one could find the right way to look at it' and suggested VV look for the `right box' for the problem. Though mysterious at the time, happily the current work is the result of understanding the insight in that stray comment.

One of the motivations of the present work is to discern whether the reconstruction is better behaved in the shallow-water regime and if reconstruction is possible over a wide range of values for the shallowness parameter. Additionally we seek to recover the bottom profile from measurements of the surface deviation $\eta$ alone. Similar to \cite{vasan2013inverse} we only assume the velocities are consistent with the shallow-water wave model and obtain the velocities as a by-product of our reconstruction. Unlike the shallow-water reconstruction method of \cite{kar2018ocean}, we only consider shallow-water wave models that are dispersive. Not only do dispersive PDEs offer some mathematical advantages over their non-dispersive counterparts in terms of smoothing and regularity, there are theoretical and experimental arguments that substantiate the need to take into account dispersion and the strong non-hydrostatic effects induced by varying bottom-boundaries \cite{dellar2005shallow,grue1992nonlinear,davies1984surface} even in the case of depth-averaged models \cite{camassa1996long}. 

\begin{figure}
\centering
\subcaptionbox{Profile 1\label{fig:profile_wavy}}[0.4\linewidth]{\includegraphics[width=0.4\textwidth]{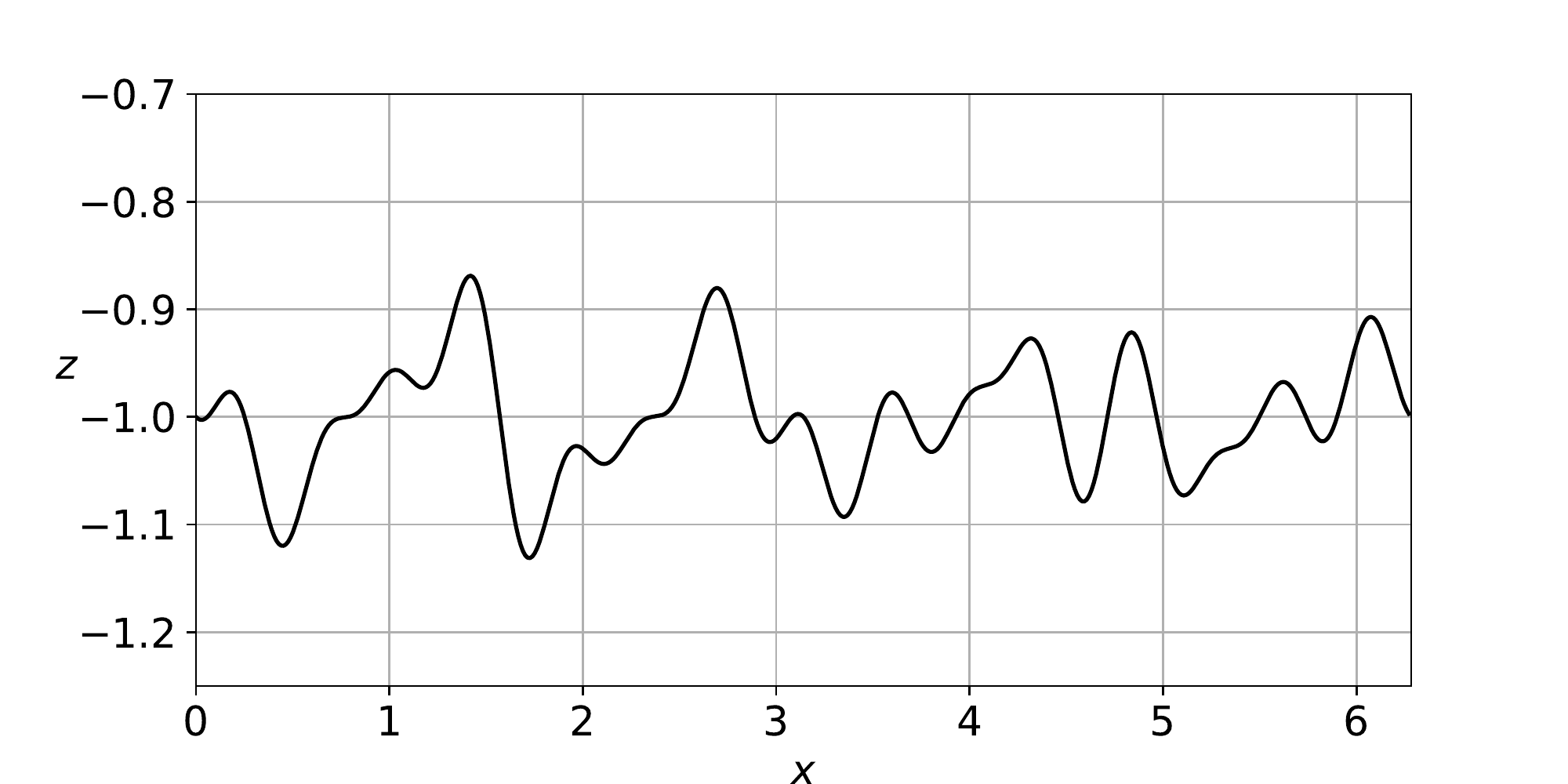}}
\subcaptionbox{Profile 2\label{fig:profile_isolated}}[0.4\linewidth]{\includegraphics[width=0.4\textwidth]{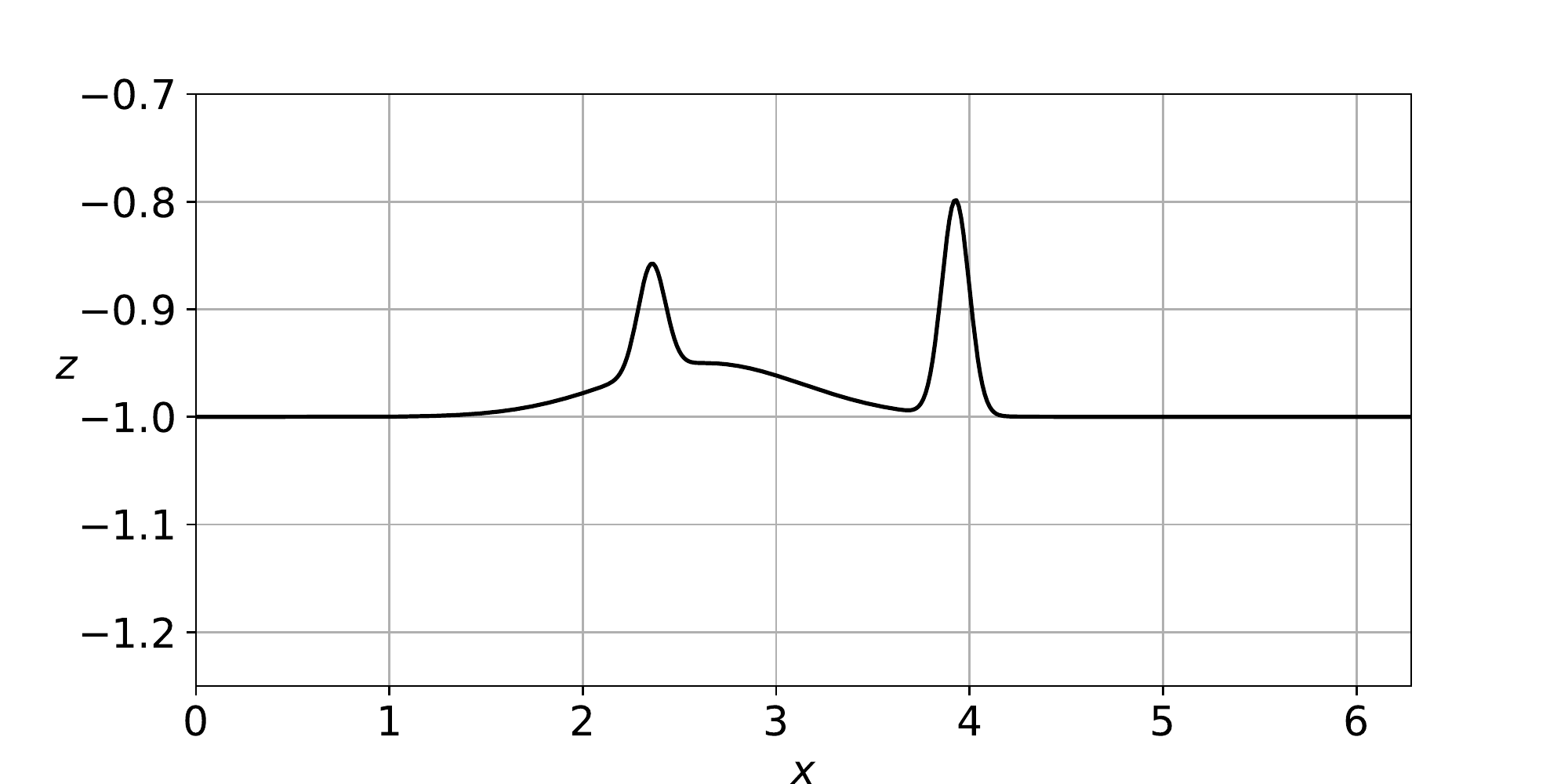}}
\caption{Two different bottom-boundary profiles. Note the scale in vertical and horizontal directions is not the same. See equations (\ref{eqn:profile1}-\ref{eqn:profile2}) for the exact form of the boundary profiles.}\label{fig:profiles}
\end{figure}

The outline of the paper is as follows. In Section \ref{sec:shallow-ww} we sketch a derivation for a class of bi-directional dispersive Boussinesq-type shallow-water wave models with varying bottom-boundaries. We focus our efforts on two such models, one of which is known in the literature and another which is, to our knowledge, novel. Following this, in Section \ref{sec:bathymetry_reconstruct} we describe the first of our inverse problems: the reconstruction of the bottom-boundary profile from surface data. In this section and the remainder of the paper, we consider the reconstruction of two specific bottom-boundary profiles representing non-isolated and isolated topographic features (see also Figure \ref{fig:profiles})
\begin{align}
\mbox{Profile 1 }:\quad \zeta = -0.12\sin(3x)\cos(2x)\sin(10x)+0.05\sin(4x),\label{eqn:profile1}
\end{align}
\begin{align}
\mbox{Profile 2 }:\quad \zeta = - 0.1e^{-100(x-x_1)^2} - &0.05e^{-2(x-x_2)^2} - 0.2e^{-100(x-x_3)^2}, \label{eqn:profile2} \\
x_1 = 3\pi/4, \: x_2 = &1.12x_1,\: x_3 = 5\pi/4. \nonumber
\end{align}
Then in Section \ref{sec:observers} we employ the observer framework \cite{luenberger1971introduction} to determine the velocity of the fluid given the surface deviation, when the bottom-boundary profile is known. The observer framework is popular in the geophysics community and is routinely used in data assimilation to recover the state of a dynamical system from partial or sparse measurements \cite{ApteA10,ApteAurouxRamaswamy2018observers,Auroux06,auroux-bonnabel11}. We then combine the two inverse problems in Section \ref{sec:simultaneous} to design our algorithm for bottom-boundary reconstruction from surface-deviation measurements alone. We do not require the mean-depth to be known and additionally we are able to reconstruct bottom profiles in deeper water than in \cite{vasan2013inverse} (though still in the shallow regime). Finally we summarise our method and describe possible extensions of the present work in Section \ref{sec:summary}.

Two key theoretical points are not discussed in the present work but nonetheless warrant attention. Firstly, one would like to know to whether the shallow-water models we consider actually model or simulate the dynamics of the full water-wave problem. For one of our models, this is already known in the literature \cite{dinvay2019well}. The other model bears a close resemblance to known Boussinesq-type equations and we conjecture the ideas of \cite{bona2005long,lannes2008justifying} may be adapted to obtain the required theorem. Secondly, one would like to know whether the recovered bottom-boundary in the shallow-water regime is close to the bottom-boundary obtained from the associated inverse problem for the full water-wave equations, such as what \cite{fontelos2017bottom} consider. Essentially this boils down to showing the approximate DNO one employs in the shallow-water model approximates the DNO of the full water-wave problem. Though we do not have rigorous statements for the models in consideration, we note results of this type are available for a variety of long-wave models \cite{lannes2013water}. We hope the success of our reconstruction algorithm spurs interest in these theoretical questions too.


\section{Shallow-water wave equations}\label{sec:shallow-ww}
Starting from the continuum description of fluid flow, the equations governing the motion of an inviscid incompressible irrotational constant-density fluid are given by 
\begin{align}\label{eqn:ww1}
	\phi_{xx} + \phi_{zz} &= 0,\quad &&-h-\zeta(x) < z < \eta(x,t),\: 0<x<L,\\
	\phi_z + \zeta_x \phi_x &=0,\quad &&z=-h-\zeta(x),\\
	\eta_t &= \phi_z - \eta_x \phi_x,\quad &&z=\eta(x,t),\\
	\phi_t + \frac 1 2 \left( \phi_x^2 + \phi_z^2\right) + g \eta &=0,\quad &&z=\eta(x,t).\label{eqn:ww2}
\end{align}
Here $\phi$ represents the velocity potential (and hence the fluid velocity is given by the gradient of $\phi$), $\eta(x,t)$ represents the free surface deviation, $\zeta(x)$ is the shape of the bottom boundary, $h$ is a typical depth, $L$ is the lateral extent of the fluid and $g$ is the acceleration due to gravity. The equations above are supplemented with periodic boundary conditions in the $x-$variable. Though we state the equations for a fluid with only one horizontal variable, the equations may just as easily be written for the more realistic scenario with two horizontal dimensions. To simplify our discussion, we restrict ourselves to the equations as specified above. 

As noted by \cite{Zakharov,CraigSulem}, equations (\ref{eqn:ww1}-\ref{eqn:ww2}) have an equivalent Hamiltonian description 
\begin{align}\label{eqn:HamiltonianWW}
	\eta_t = \frac{\delta H}{\delta q},\quad q_t = -\frac{\delta H}{\delta \eta},
\end{align}
where the Hamiltonian is given, in terms of the surface deviation $\eta$ and the Dirichlet-trace of the velocity potential at the surface $q(x,t)=\phi(x,\eta,t)$, by
\[H = \frac 1 2 \int \left(q G(\eta,\zeta)q + g\eta^2\right) \: dx.\]
Here $G(\eta,\zeta)$ is the Dirichlet-Neumann operator (DNO) that maps the Dirichlet condition of the following boundary-value problem
\begin{alignat}{3}
	\psi_{xx} + \psi_{zz} &= 0, &\quad  -&h-\zeta(x) < z < \eta(x,t),\: 0<x<L,\\
	\psi_z + \zeta_x \psi_x &=0, &\quad z&=-h-\zeta(x),\\
	\psi &= q(x,t), &\quad z&=\eta(x,t),
\end{alignat}
to the associated Neumann condition at the surface $z=\eta(x,t)$. Thus \[G(\eta,\zeta)q = \psi_z - \eta_x \psi_x,\quad z=\eta(x,t).\]
For the problem posed on the whole line, for any $k_0>1/2$, if $\eta,\zeta$ are elements of the (real-valued) Hilbert space $H^{k_0+1}(\mathbb{R})$, at each time $t$, such that there exists a constant $h_0>0,\: \eta(x,t)+h+\zeta(x)\geq h_0$, then the Dirichlet-Neumann operator is a bounded linear operator defined as $G(\eta,\zeta): \dot{H}^{3/2}(\mathbb{R})\to H^{1/2}(\mathbb{R})$ \cite{lannes2013water} where $\dot{H}^{3/2}(\mathbb{R})$ is a Beppo-Levi space \cite{deny1954espaces}. Moreover $G(\eta,\zeta)$ has a self-adjoint realisation on $L^2(\mathbb{R})$ with domain $H^1(\mathbb{R})$ (see Section $3.1$ of \cite{lannes2013water}). On the other hand, if $\eta(x,t) = \tilde \epsilon f(x,t)$ and $\zeta(x) = \tilde\delta g(x)$ where $f,g$ are real-analytic functions of $x$ (for every $t$), then $G(\eta,\zeta)$ is analytic in $\tilde\epsilon,\tilde\delta$ \emph{i.e.} $G$ has a convergent Taylor-series operator-expansion when $q$ is also real-analytic in $x$ (for every $t$) \cite{nichollsTaber2008joint}.

We emphasise two points. Firstly, equations (\ref{eqn:HamiltonianWW}) imply the dynamically relevant quantities are those defined on the boundary: the shape of the free surface $\eta$ and the value of the potential at the surface $q(x,t)=\phi(x,\eta,t)$. Secondly, although the equations look like a standard partial differential equation, we note that $G(\eta,\zeta)$ is \emph{not} a local operator.

We now describe the formal procedure to derive models whose solutions approximate the full fluid motion. The overall methodology we follow was introduced in \cite{benjamin1984impulse} and employed in \cite{broer1974hamiltonian,broer1975approximate,broer1976stable} to deduce model water-wave equations. Model equations using approximate Hamiltonians were also employed by \cite{craig1994hamiltonian} though we consider varying bottom topography. Essentially one replaces the DNO $G(\eta,\zeta)$ by a simpler operator $G^M(\eta,\zeta)$ (which we refer to as a model DNO). Note we do \emph{not} obtain effective equations via a homogenisation theory as described in \cite{Craigetal,craig2009water}. 

One way of constructing approximations to the DNO is due to \cite{afm} where they characterise the DNO in terms of the following equations
\begin{align}\label{eqn:AFM_traditional}
	\int e^{ikx}\left\{ iG(\eta,\zeta)q\:\mathrm{cosh}(k(\eta+h)) + q_x\:\mathrm{sinh}(k(\eta+h)) + Q_x\:\mathrm{sinh}(k\zeta) \right\}dx = 0,\\
	\int e^{ikx}\left\{ iG(\eta,\zeta)q\:\mathrm{sinh}(k(\eta+h)) + q_x\:\mathrm{cosh}(k(\eta+h)) - Q_x\:\mathrm{cosh}(k\zeta) \right\}dx = 0.
\end{align}
Here $Q_x$ is the velocity tangential to the fluid domain at the bottom $z=-h-\zeta$: $Q_x = \phi_x - \zeta_x\phi_z$ evaluated at the bottom boundary. We assume that $\|\zeta\|_{\infty}/h$ is small. Given $q_x$, one solves the above equations simultaneously for both $G(\eta,\zeta)q$ and $Q_x$. Note also that we do not assume the mean of $\zeta$ is zero. In this case, $h$ can be thought of as an estimate for the bottom topography rather than the actual mean bottom surface. 

Within the shallow-water regime, the aspect ratio $h/L$ is assumed to be a small parameter. Suppose further that \emph{both} 
\[\frac{\|\eta\|_\infty}{h},\: \frac{\|\zeta\|_\infty}{h} \sim \left(\frac{h}{L}\right)^2.\]
Then we have the following $\mathcal{O}(\,(h/L)^4)$ accurate expansion for the hyperbolic functions
\begin{align}
	\mathrm{cosh}(k(\eta+ h)) &= 1 + \frac{(kh)^2}{2}  + \ldots,\\
	\mathrm{sinh}(k(\eta+h)) &= kh + k \eta + \frac{ k^3h^3}{6} + \ldots,\\
	\mathrm{cosh}(k \zeta) &= 1 + \ldots,\\
	\mathrm{sinh}(k\zeta) &= k \zeta + \ldots.
\end{align}
Using these expansions in the global relations (\ref{eqn:AFM_traditional}) we obtain 
\begin{align}
\int e^{ikx} \left\{ \left(1 + \frac{(kh)^2}{2}\right) iG(\eta,\zeta)q + q_x \left( kh + k \eta + \frac{ k^3h^3}{6} \right) + Q_x k \zeta  \right\}dx = \mbox{h.o.t.'s\:},\\
\int e^{ikx} \left\{ \left( kh + k \eta + \frac{ k^3h^3}{6} \right) iG(\eta,\zeta)q + q_x \left(1 + \frac{(kh)^2}{2}\right)  - Q_x  \right\}dx = \mbox{h.o.t.'s\:}.
\end{align}
The higher order terms in the above expressions involve terms of $\mathcal{O}(\,(h/L)^4)$ and higher. To obtain an expression for the DNO to similar order, we need only obtain an expression for $Q_x$ which is $\mathcal{O}(1)$ accurate. This is readily obtained from the second of the above equations: $Q_x = q_x$ which is consistent with our shallow-water approximation. This leads to the following expression for a model DNO with variable bottom-boundary 
\begin{align}
	\left( 1 - \frac{h^2}{ 2} \partial_x^2\right) G^M(\eta,\zeta)q = -hq_{xx} + \frac{h^3}{6}q_{xxxx} - \partial_x[ (\eta+\zeta) q_x].
\end{align}
Note however that the model DNO given above is not formally self-adjoint. Indeed on solving for $G^M$ explicitly we have
\begin{align}
 G^M(\eta,\zeta)q = -\left( 1 - \frac{h^2}{ 2} \partial_x^2\right)^{-1}\left(h - \frac{h^3}{6}\partial_x^2\right)q_{xx} - \left( 1 - \frac{h^2}{ 2} \partial_x^2\right)^{-1}\partial_x[ (\eta+\zeta) q_x],
\end{align}
which, due to the last term on the right-hand side, is not a self-adjoint operator. To the same level of asymptotic accuracy one could consider the alternate operator
\begin{align}
	 G^M(\eta,\zeta)q = -\left( 1 - \frac{h^2}{ 2} \partial_x^2\right)^{-1}\left(h - \frac{h^3}{6}\partial_x^2\right)q_{xx}  - \partial_x[ (\eta+\zeta) q_x],
\end{align}
which \emph{is} a formally self-adjoint operator acting on $q$. With the model DNO $G^M$ defined as above, we construct a Hamiltonian system
\begin{align}
\eta_t = \frac{\delta H^M}{\delta q} \quad q_t = -\frac{\delta H^M}{\delta \eta},\quad \quad
	H^M = \frac 1  2 \int \left(qG^M(\eta,\zeta)q + \eta^2 \right)dx. \label{eqn:modelHamil}
\end{align}
This gives rise to an equation which is essentially the same Boussinesq-type equation derived by \cite{afm}, albeit with a non-trivial bottom boundary. However,  we prefer to work with a regularised version of this Boussinesq equation which has the same level of formal asymptotic accuracy
\begin{align}\label{eqn:omegaPa}
    \eta_t &= \omega^2(-i\partial_x)q - \mathcal{P}(-i\partial_x) \partial_x\left[ (\eta+ \zeta)\mathcal{P}(-i\partial_x)q_x\right],\\ \label{eqn:omegaPb}
	q_t &= -\eta - \frac 1 2 \left(\mathcal{P}(-i\partial_x)q_x \right)^2,
\end{align}
where 
\begin{align}\label{eqn:regBouss}
\omega^2(k) = \frac{1+\frac{(\mu k)^2}{6}}{1+\frac{(\mu k)^2}{2}}k^2,\quad \mathcal{P}(k) = \frac{1}{1+\frac{(\mu k)^2}{2}},\quad \mu = \frac{2\pi h}{L},
\end{align}
and we have also nondimensionalised the equations with a horizontal length-scale $L/(2\pi)$; a scale $h$ for both $\zeta,\eta$; a scale $\sqrt{gh}$ for $q_x$ and chosen the longer time-scale $L/(2\pi \sqrt{gh})$. The factors of $2\pi$ are for convenience so that our wavenumbers are integers. With this scaling, the only non-dimensional parameter is the aspect ratio $\mu$. Equations (\ref{eqn:omegaPa}-\ref{eqn:omegaPb}) are Hamiltonian with canonical Poisson structure with the Hamiltonian given by
\begin{align}
	H^{\omega,\mathcal P} = \frac  1 2\int\left( q\omega^2(-i\partial_x)q + (\eta+\zeta)(\mathcal P q_x)^2 + \eta^2\right)\: dx.
\end{align}
Note if we make a further substitution $\omega^2 = k^2 + \mathcal{O}(\mu^2)$ and $\mathcal P =  1 + \mathcal{O}(\mu^2)$ we obtain a different Hamiltonian system: the hydrostatic shallow-water equations in one horizontal dimension. Various choices for the pseudo-differential operators $\omega^2,\mathcal P$ can lead to the non-trivial bottom-boundary versions of the ASMP model \cite{aceves2013numerical} or a Hamiltonian version of the Hur-Pandey model \cite{hur2019modulational}. All these models are different types of bidirectional Whitham equations. For a comparison of different such models see \cite{dinvay2019comparative}. 

For the remainder of this paper, we concern ourselves with only two specific shallow-water wave models, both given by equations (\ref{eqn:omegaPa}-\ref{eqn:omegaPb}). The first we refer to as regularised Boussinesq and is given by the choice (\ref{eqn:regBouss}). The second model equation is given by the choice
\begin{align}\label{eqn:regBoussWhitham}
\omega^2(k) = k\:\frac{\mathrm{tanh}(\mu k)}{\mu},\quad \mathcal{P}(k) =\frac{\mathrm{tanh}(\mu k)}{\mu k},\quad \mu = \frac{2\pi h}{L},
\end{align}
which leads to the equation considered by \cite{dinvay2019well} which we refer to as regularised Boussinesq-Whitham (which is also the Hamiltonian version of the Hur-Pandey model). Note from here on, for the sake of notational convenience we will suppress the argument of $\omega^2$ and $\mathcal P$, with the understanding that they are suitably interpreted either as an operator acting on functions of the real variable $x$ or a multiplier acting on Fourier transforms of functions.

{\remark{The equations considered in \cite{carter2020fully} are essentially a generalised version of (\ref{eqn:omegaPa}-\ref{eqn:omegaPb}) which include a boundary operator in the place of $\zeta$. The shallow-water wave models described in the current work involve shallow-water approximations to the boundary operator employed in \cite{carter2020fully}. }}
{\remark{To make sense of the pseudo-differential operators we restrict ourselves only to the problem with periodic boundary conditions. In a sense, this assumption was adopted when we claimed the Hamiltonian formulation of the full water-wave problem. A suitable phase space needs to be chosen before taking a variational derivative. The extension of the current work and equations such as (\ref{eqn:omegaPa}-\ref{eqn:omegaPb}) to non-periodic domains is interesting but will be left for future works.}}
{\remark{The regularised Boussinesq-Whitham equations with trivial bottom-boundary was considered in \cite{dinvay2019well} where they showed local and global wellposedness. Their equations are slightly different from those given here. Indeed their equations are written in terms of a new velocity variable. It suffices for our purposes, and makes our computer code more modular, to consider equations in terms of velocity potential since only $\omega$ and $\mathcal P$ need to be defined.}}
{\remark{Both regularised Boussinesq and regularised Boussinesq-Whitham have equations where the nonlinearity is a bounded operator (on some suitable function space such as $L^2([0,2\pi])\times L^2([0,2\pi])$). Indeed regularised Boussinesq has a smoothing nonlinearity. For this reason we expect regularised Boussinesq to possess a local wellposedness theory for sufficiently smooth initial data.}}
{\remark{When employing a canonical Poisson bracket to define Hamiltonian evolution equations, as was done in the above passage, only the symmetric part of the model $DNO$ appears in the equations of motion. Thus if we require both the Hamiltonian and the dynamics of the model to approximate the full water-wave system, we impose the model DNO $G^M$ to be a symmetric operator on $q$. This property has the further advantage that any system of the form (\ref{eqn:modelHamil}) where $G^M$ is symmetric, automatically conserves the momentum $I = \int q_x\eta$ when $\zeta=0$. This is evident when one computes the (canonical) Poisson bracket of the Hamiltonian $H^M$ and the momentum $I$. Likewise additional properties of $G^M$ imply further conserved quantities: if the range of $G^M$ is orthogonal to constant functions, then $\int \eta$ is conserved; if $G^M$ is symmetric and $G^M \:x=-\eta_x$ in a distributional sense, then the centre of mass $\int x\eta - t \int q_x\eta$ is conserved and so on. We do not pursue this point any further since the additional conserved quantities do not directly have a bearing on our problem.}}


\section{Bathymetry using surface data}\label{sec:bathymetry_reconstruct}
In the bottom-profile reconstruction approach of \cite{fontelos2017bottom}, the authors determine the profile $\zeta$ as the minimiser of the functional \[F(\zeta) = \int_0^{2\pi} (\eta_t - G(\eta,\zeta)q )^2dx, \]
where $\eta,q,\eta_t$ are given functions. They show that for the DNO of the water-wave problem, a minimiser exists and is in fact unique, when the data $\eta,\eta_t,q$ are given from a solution to the water-wave equations (\ref{eqn:ww1}-\ref{eqn:ww2}). This approach has an obvious reformulation for shallow-water wave models which employ a model DNO: one replaces the $G$ in the above equation by a suitable approximate DNO $G^M$. For model equations of the form (\ref{eqn:omegaPa}-\ref{eqn:omegaPb}) this leads to
\begin{align}\label{eqn:bathymetry_min}
\zeta = \argmin_{\zeta^*} \int(\eta_t - \omega^2 q + \mathcal P \partial_x((\eta+\zeta^*)\mathcal P q_x)\:)^2dx.
\end{align}
The same functional is obtained if one adapts the reconstruction algorithm of \cite{vasan2013inverse} to the shallow-water regime. For details, see Appendix \ref{app:VD_shallow}. Note that the operator $\mathcal P$, for either choice (\ref{eqn:regBouss}) or (\ref{eqn:regBoussWhitham}), maps $L^2([0,2\pi])$ functions to at least $H^1([0,2\pi])$. If we further assume $\eta_t,\eta\in L^2([0,2\pi])$ and $q_x\in H^1([0,2\pi])$, then a reasonable function space for the above minimisation problem is $\zeta^*\in L^2([0,2\pi])$. This suffices for both regularised Boussinesq and regularised Boussinesq-Whitham. Before we proceed, we emphasise the notation used in the rest of the paper. The expression $\mathcal Pq_x$ appears repeatedly in the following. It will always denote a {\textit{function}}. When this expression appears in the definition of an operator, it denotes multiplication by the function $\mathcal Pq_x$.

We now suppose that $\eta,\eta_t,q_x$ are known functions. Then to find the minimiser of (\ref{eqn:bathymetry_min}) we write the associated Euler-Lagrange equation
\begin{align}\label{eqn:bathymetry_EL}
	(\mathcal Pq_x)\: \mathcal P^2\partial_x^2 \left( (\mathcal Pq_x) \: \zeta^*\right) = -(\mathcal Pq_x)\: \mathcal P\partial_x 
	\Big( \eta_t - \omega^2 q + \mathcal P\partial_x \left( (\mathcal Pq_x)\: \eta \right) \Big).
\end{align}
Evidently $\zeta^*=\zeta$ is a minimiser of (\ref{eqn:bathymetry_min}). From equations (\ref{eqn:omegaPa}) and (\ref{eqn:bathymetry_EL}) we find
\begin{align}
	(\mathcal Pq_x)\: \mathcal P^2\partial_x^2 \left( (\mathcal Pq_x) \: \zeta^*\right) = (\mathcal Pq_x)\: \mathcal P^2\partial_x^2 \left( (\mathcal Pq_x )\: \zeta\right),
\end{align}
thus $\zeta$ (the true bottom boundary) is a solution to (\ref{eqn:bathymetry_EL}). To solve equation (\ref{eqn:bathymetry_EL}) for the bottom-profile $\zeta^*$ we need to invert the operator \[B:f \to (\mathcal Pq_x) \mathcal P^2\partial_x^2\big( (\mathcal Pq_x)\:f\big).\]
This operator is formed by compositions of two operators: multiplication by the function $\mathcal Pq_x$ and the operator $\mathcal{P}^2\partial_x^2$. Note the latter corresponds to the Fourier multiplier $-k^2\mathcal{P}(k)^2$. Recall the definition of $\mathcal P$ depends on the particular shallow-water model we use. As mentioned above, the operator $\mathcal P$, for either model, maps $L^2([0,2\pi])$ functions to at least $H^1([0,2\pi])$. Thus from a standard Sobolev embedding the function $\mathcal Pq_x$ is continuous in $x$ when $q_x\in L^2([0,2\pi])$. Thus multiplication by $\mathcal  Pq_x$ is a bounded operator on $L^2([0,2\pi])$. Moreover, for either model, $k^2\mathcal{P}(k)^2$ is bounded. It then follows that the operator $B$ is a bounded linear operator on $L^2([0,2\pi])$. It is readily verified $B$ is self-adjoint. Next, $B$ has spectrum contained in the negative real axis since
\[\int f Bf\:dx = \int f \: (\mathcal Pq_x)\: \mathcal P^2\partial_x^2 \left( (\mathcal Pq_x) \: f \right) \:dx = - \int \left(\mathcal P\partial_x \big( (\mathcal Pq_x)\:f \big)\right)^2\:dx. \]
In the case of regularised Boussinesq, $k^2\mathcal{P}(k)^2\sim k^{-2}$ for large $|k|$. This implies $\mathcal{P}^2\partial_x^2$ and hence $B$, maps $L^2([0,2\pi])$ to $H^2([0,2\pi])$. This implies $B$ is in fact a compact operator as $H^2([0,2\pi])$ is compactly embedded in $L^2([0,2\pi])$ \cite{evans1998partial}. On the other hand, the operator $B$ is not compact for regularised Boussinesq-Whitham. Indeed we have the following theorem.

{\theorem{Suppose $q_x\in L^2([0,2\pi])$ and $\mathcal Pq_x$ is not identically zero, then for regularised Boussinesq-Whitham, $B$ is not a compact operator.
\begin{proof}
	The operator $\mathcal{P}^2\partial_x^2$ has the Fourier symbol $-\tanh^2(\mu k)/\mu^2$ which is bounded but does not decay for large $|k|$. As a result, the operator $\mathcal{P}^2\partial_x^2$ is not compact. Indeed consider the orthonormal sequence in $L^2([0,2\pi])$, $v_n = \sin(nx)/\sqrt{\pi},\: n=1,2,\ldots$ with $\|v_n\|_2=1$. Then 
	\[\|\mathcal{P}^2\partial_x^2v_n-\mathcal{P}^2\partial_x^2v_m\|_2^2 = (\tanh(\mu n)/\mu)^4 + (\tanh(\mu m)/\mu)^4 \geq 2 (\tanh(\mu)/\mu)^4.\]
	We conclude $v_n$ is a bounded sequence that satisfies 
	\[\|\mathcal{P}^2\partial_x^2(v_n-v_m)\|_2\geq \sqrt 2 (\tanh(\mu)/\mu)^2\mbox{ for every }n\neq m.\] 
	This implies no subsequence of $\mathcal P^2\partial_x^2v_n$ converges. Thus the operator $\mathcal{P}^2\partial_x^2$ is not compact. 

	We now show that $B$ is not compact. The function $\mathcal Pq_x$ is continuous when $q_x\in L^2([0,2\pi])$ for regularised Boussinesq-Whitham. Let $S$ be a non-empty connected component of the set $\{x\in[0,2\pi]: |\mathcal Pq_x|\geq c\}$ for some $c>0$. Since $\mathcal Pq_x$ is not identically zero, there exists such a $c$.

	Note $\mathcal Pq_x$ restricted to $S$ is in $L^2(S)$. Let $X\subset L^2(S)$ be the orthogonal complement to the span of $\mathcal Pq_x$ and take an orthonormal sequence $w_n$ in $X$. The sequence is obtained from an orthonormal sequence in $L^2(S)$, projected onto $X$ followed by the Gram-Schmidt process and then considering only non-trivial $w_n$. Define $v_n\in L^2[0,2\pi]$ by $v_n= w_n \chi_S$ where $\chi_S$ is the indicator function on $S$. Then $v_n$ is a bounded sequence in $L^2([0,2\pi])$ whose image under $B$ satisfies 
	\begin{align*}
	\|Bv_n-Bv_m\|_2^2 &= \int_0^{2\pi}\left( (\mathcal Pq_x)\: \mathcal P^2\partial_x^2 \big( (\mathcal Pq_x) \: (v_n-v_m)\big) \right)^2dx,\\
	&\geq \frac{c^2}{\mu^4} \sum_{k=-\infty}^\infty \left(\tanh^4(\mu k)\: |\mathcal F\left[(\mathcal Pq_x) \: (v_n-v_m)\right]_k|^2\right),\\ 
	&\geq 2 c^4 (\tanh(\mu)/\mu)^4,\quad n\neq m,
	\end{align*}
	where in the second line we use Parseval's identity. To obtain the third line, we use Parseval's in the reverse direction for which the orthogonality condition $\int_S  (\mathcal Pq_x)w_n =0$ is crucial, as the zero index does not contribute to the summation. Once again it follows that no subsequence of $Bv_n$ converges and hence $B$ is not a compact operator.	
\end{proof}
}}
{\remark For regularised Boussinesq-Whitham, the function $\mathcal Pq_x$ is identically  zero if and only if $q_x=0$, since $\mathcal P(k) = \tanh(\mu k)/(\mu k) >0$ for all $k$. \\
}

For either shallow-water model, we have the following result regarding zero eigenvalues of the operator $B$.
{\theorem{Suppose $q\in H^2([0,2\pi])$ and the function $\mathcal Pq_x$ only vanishes on a set of measure zero, then the operator \[(\mathcal Pq_x)\: \mathcal P^2\partial_x^2 \big( (\mathcal Pq_x) \: \cdot\big)\] has no zero eigenvalue for $\mathcal P$ given by either (\ref{eqn:regBouss}) or (\ref{eqn:regBoussWhitham}).
\begin{proof}
	Suppose there is a zero eigenvalue. Then 
		\[(\mathcal Pq_x)\: \mathcal P^2\partial_x^2 \big( (\mathcal Pq_x) \: f \big) = 0,\]
	for some $f\in L^2([0,2\pi])$. But then we have 
	\[ \mathcal P^2\partial_x^2 \big( (\mathcal Pq_x) \: f \big) = 0 \Rightarrow (\mathcal Pq_x) f = C,\]
	for some constant $C$, which follows from the fact that $\mathcal P^2\partial_x^2$ in Fourier space is the multiplier $-k^2\mathcal P(k)^2$ and $\mathcal P$ is a positive operator (see definition in \ref{eqn:regBouss} or \ref{eqn:regBoussWhitham}). 

	Suppose the constant $C=0$. Since $\mathcal Pq_x$ vanishes on a set of measure zero, we conclude $f=0$ almost everywhere. But then $f$ cannot be an eigenfunction and hence there is no zero eigenvalue. Thus we assume $C\neq 0$.

	For $q\in H^2([0,2\pi])$, $\mathcal P q_x$ is given by $ik \mathcal P \hat q$ where $\hat q$ is the Fourier transform of $q$. Evidently the function $\mathcal P q_x$ has zero average. Moreover since $\mathcal P\partial_x$ is a bounded operator, $\mathcal Pq_x$ is also in $H^2([0,2\pi])$ and from a Sobolev embedding theorem, $\mathcal Pq_x\in C^{1,1/2}([0,2\pi])$. Since $\mathcal Pq_x$ is a continuous function with zero average, it must vanish somewhere. Since $\mathcal P q_x$ is periodic, we assume without loss of generality $\mathcal Pq_x=0$ at $x=0$. Then from Taylor's theorem in Lagrange form for a neighbourhood of $x=0$ we have
	\[\mathcal Pq_x = x \left.\frac{d\mathcal Pq_x}{dx}\right|_{x=s} = x \left.\frac{d\mathcal Pq_x}{dx}\right|_{x=0} + x \left( \left.\frac{d\mathcal Pq_x}{dx}\right|_{x=s} - \left.\frac{d\mathcal Pq_x}{dx}\right|_{x=0}  \right)
	,\]
	for some $s$ between $0$ and $x$, and $|x|$ sufficiently small. As the derivative of $\mathcal Pq_x$ is a $1/2$-H\"older function, this implies
	\[|\mathcal P q_x| \leq |x| m_1 + m_2 |x|^{3/2} \Rightarrow |\mathcal P q_x| \leq |x|(m_1 + m_2),\mbox{ for }|x|<1, \]
	for some constants $m_1,m_2$. Note, $m_1$ and $m_2$ cannot both be zero as that would imply $\mathcal Pq_x$ vanishes in an open set. This implies the eigenfunction $f$ is given by 
	\[f = \frac{C}{\mathcal P q_x} \Rightarrow |f| \geq \frac{|C|}{|x|(m_1+m_2)},\mbox{ for }|x|\mbox{ sufficiently small}.\]
	
	But then $f\notin L^2([0,2\pi])$ and hence zero cannot be an eigenvalue.
\end{proof}
}}
Despite the uniqueness result mentioned above, we expect the operator on the left-hand side (\ref{eqn:bathymetry_EL}) to be ill-conditioned. For regularised Boussinesq, where $\mathcal P$ is given by (\ref{eqn:regBouss}), as the operator $B$ is compact, zero is a limit point of the (negative) eigenvalues. Additionally, for both shallow-water models, the zeros of $\mathcal Pq_x$ can lead to further ill-conditioned behaviour. Indeed consider a sequence of smooth functions  $f_n(x)$ with compact support and which approximate a Dirac delta distribution located at a zero of $\mathcal Pq_x$. Then for any $\zeta\in L^2([0,2\pi])$ \[\lim_{n\to \infty} \int f_n(\mathcal Pq_x) \mathcal P^2\partial_x^2 \big( (\mathcal Pq_x)\: \zeta\big) \:dx = 0.\]
This holds for both regularised Boussinesq and regularised Boussinesq-Whitham and is unavoidable. Indeed for $q\in H^1([0,2\pi])$, $\mathcal Pq_x$ is continuous and has zero average, which implies $\mathcal Pq_x$ must be zero for some $x\in[0,2\pi]$. Thus the main difficulty in our approach to bathymetry is taming the ill-conditioned nature of the operator $B$. This ill-conditioned nature of $B$ is directly related to the ill-posed nature of bottom-boundary detection in the context of the water-wave problem (\ref{eqn:ww1}-\ref{eqn:ww2}).

{\remark{}
So far we have only discussed reconstruction with relation to the dispersive models of (\ref{eqn:regBouss}) and (\ref{eqn:regBoussWhitham}). We now briefly consider the hyperbolic model corresponding to the choice $\omega^2=k^2$ and $\mathcal P =1$. Now $\eta,q$ evolve according to the Saint-Venant equations. The corresponding minimisation problem for the bottom-profile $\zeta$ leads to the following Euler-Lagrange equation 
\[q_x \partial_x^2( q_x \zeta^*) = -q_x \partial_x (\eta_t + q_{xx} + \partial_x(\eta q_x)\: ).\]
For sufficiently smooth $q_x$ (at least twice differentiable), the operator on the left-hand side in the above equation is in fact a Sturm-Liouville operator \[q_x\partial_x^2 (q_x \zeta^*) = \partial_x(q_x^2 \partial_x \zeta^*) + q_x q_{xxx} \zeta^*.\]
We note the same issue regarding the zeros of the coefficient $q_x$ plagues this operator, leading to a singular SL problem. However the situation is arguably worse since the forward model for $q_x$ and $\eta$ is now a nonlinear hyperbolic PDE and likely forms shocks. Consequently, the coefficients of the operator above may be discontinuous thereby precluding any traditional SL theory. If we instead formally cancel the respective terms on either side of the above equation we obtain
\[-\partial_x q_x \zeta^* = \eta_t + \partial_x (\eta+1) q_x.\]
Further assuming that $\eta_t$ has zero average, we formally solve for $\zeta^*$ by having to divide by $q_x$. This seems a highly suspect way of obtaining the bottom boundary that is unlikely to be robust to noise. For these reasons we do not consider the hyperbolic model any further.}

\vspace{10pt}
The upshot of the discussion in the previous paragraphs is the following: to solve (\ref{eqn:bathymetry_EL}) for the bottom boundary, one needs some form of regularisation so that the eigenvalues of any finite dimensional truncation of the operator are sufficiently far away from the origin. After exploring many different strategies we concluded the most straightforward and physically meaningful approach was to recover the bottom boundary from the following minimisation problem
\begin{align}\label{eqn:bathymetry_min_alltime}
\zeta = \argmin_{\zeta^*}\sum_{j=1}^M \int(\eta_t^{(j)} - \omega^2 q^{(j)} + \mathcal P \partial_x((\eta^{(j)}+\zeta^*)\mathcal P q_x^{(j)})\:)^2dx,
\end{align}
where the superscript $j$ indicates data ($\eta,\eta_t,q_x$) obtained at time $t_j$. Thus we demand the bottom boundary $\zeta$ minimises the functional of \cite{fontelos2017bottom} at multiple time instances simultaneously. This was found particularly helpful when the zeros of $\mathcal Pq_x$ do not remain fixed in space for a given solution to the shallow-water wave model equations. We note this is not always the case, but that it is fairly straightforward to generate initial conditions that leads to such preferred solutions. Indeed an initial condition inspired by a Stokes expansion 
\begin{align}
\label{init:conditions}
q(x,0) = A \sin(x) + A/5\sin(2x),\: \eta(x,0) = A \cos(x) + A/5\cos(2x) - 0.1,
\end{align} often leads to solutions that look like travelling waves even over non-trivial bottom-boundaries. We found this initial condition was sufficient for the purpose of preventing any zero of $\mathcal Pq_x$ remain fixed in space for all time. The associated Euler-Lagrange equation for (\ref{eqn:bathymetry_min_alltime}) is similar to (\ref{eqn:bathymetry_EL}) but includes a summation over times $t_j$ on either side
\begin{align}\label{eqn:bathymetry_EL:all_time}
	\sum_{j=1}^M\left[(\mathcal Pq_x^{(j)})\: \mathcal P^2\partial_x^2 \left( (\mathcal Pq_x^{(j)}) \: \zeta^*\right)\right] = -\sum_{j=1}^M\left[\mathcal (Pq_x^{(j)})\: \mathcal P\partial_x \left(\eta_t^{(j)} - \omega^2 q^{(j)} + \mathcal P\partial_x \left( (\mathcal Pq_x^{(j)})\: \eta^{(j)} \right)\: \right)\:\right].
\end{align}

\subsection{Numerical experiments}\label{sec:bathy_reconstruct:examples}

\begin{figure}
\centering
\subcaptionbox{Regularised Boussinesq \label{fig:eigvals:bouss}(\ref{eqn:regBouss})}[0.45\linewidth]{\includegraphics[width=0.4\textwidth]{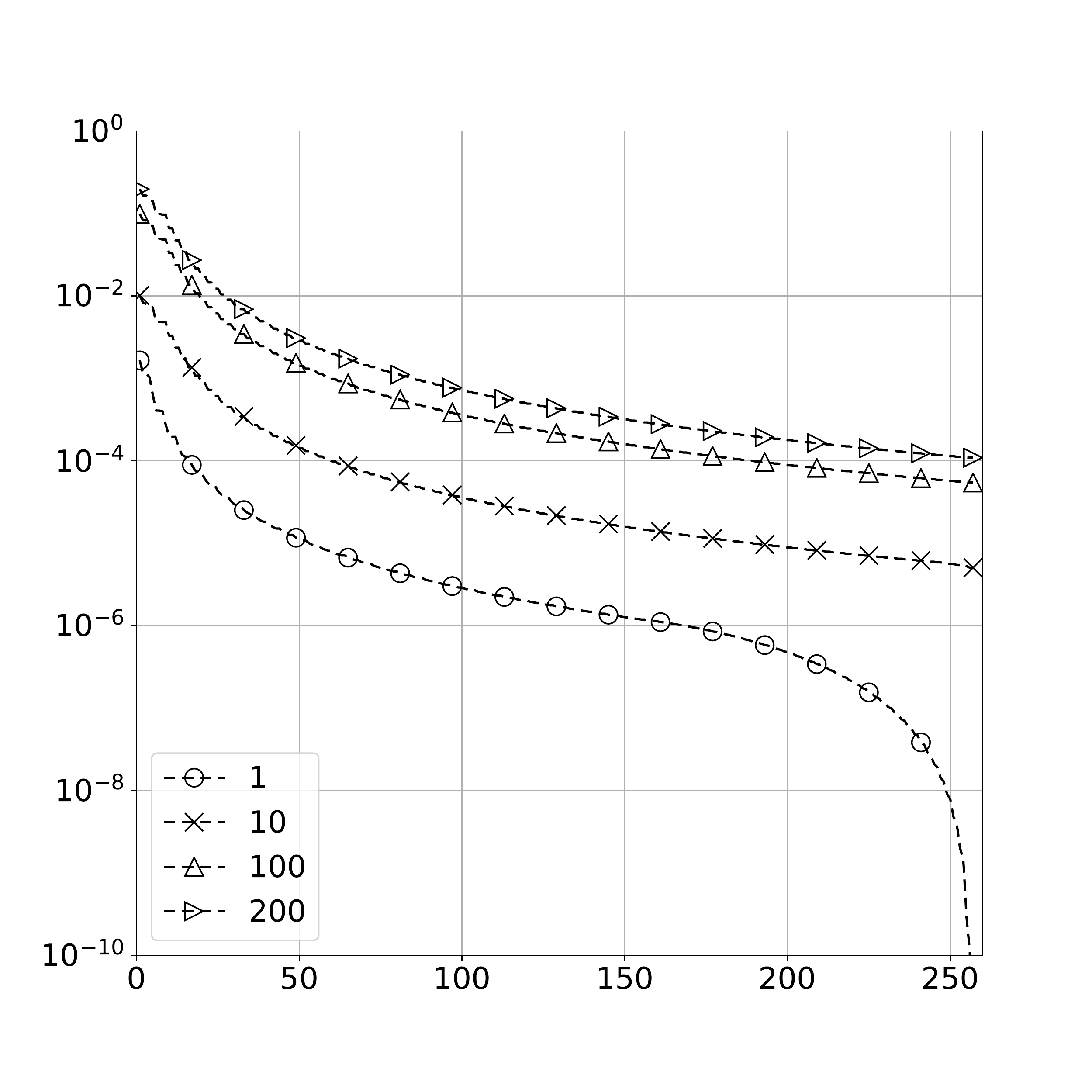}}
\subcaptionbox{Regularised Boussinesq-Whitham \label{fig:eigvals:bousswhitham}(\ref{eqn:regBoussWhitham})}[0.45\linewidth]{\includegraphics[width=0.4\textwidth]{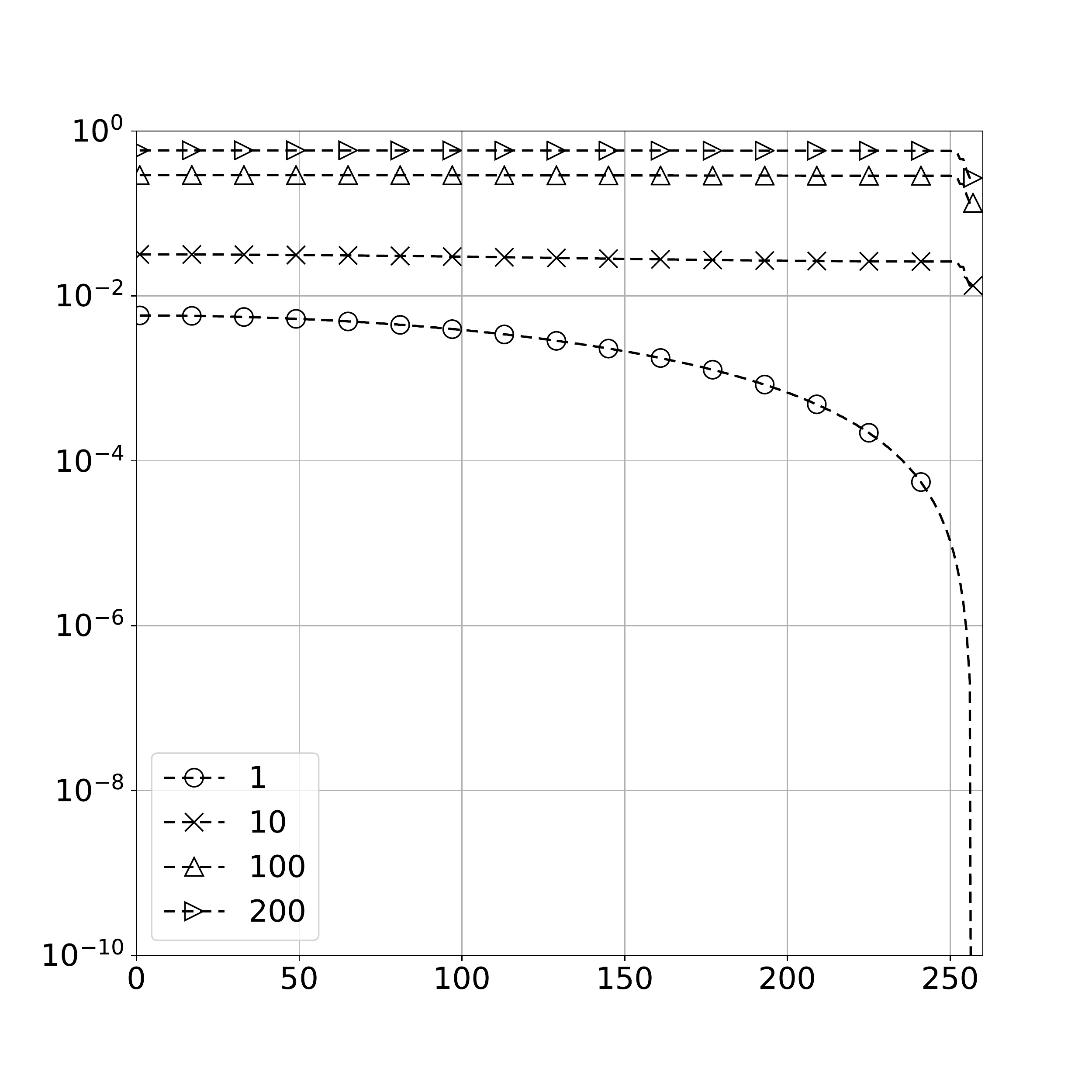}}
\caption{Absolute value of the (sorted) eigenvalues of the operator on the left-hand side of (\ref{eqn:bathymetry_EL:all_time}) for $\mu=1$ in either regularised Boussinesq or regularised Boussinesq-Whitham. Here $q^{(j)}=0.1\sin(x-t_j)$ where $x\in[0,2\pi]$ and $t_j=2\pi j/M,\: j=1,2,\ldots M$. The horizontal axis is the index of the sorted eigenvalue. The different curves correspond to different values of $M$ as indicated in the legend. The markers on the different curves are placed to distinguish the curves and do not indicate the number of data-points in the figure. The number of grid points in $x$ was $256$ in all cases for both model equations. Hence $256$ eigenvalues are shown in each curve. }\label{fig:eigvals}
\end{figure}

\begin{figure}
\centering
\subcaptionbox{Profile \ref{fig:profile_wavy} }[0.475\linewidth]{\includegraphics[width=0.475\textwidth]{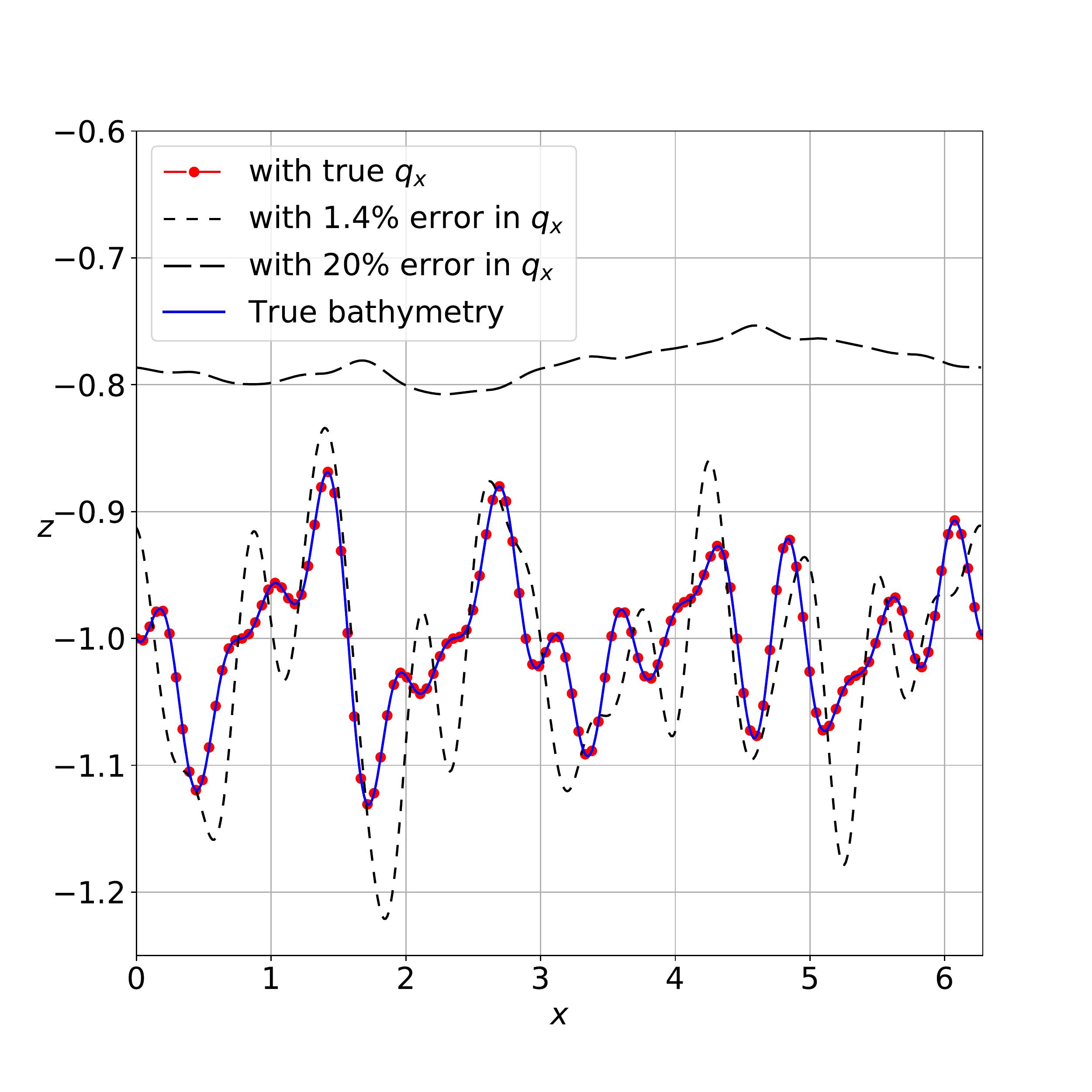}}
\subcaptionbox{Profile \ref{fig:profile_isolated}}[0.475\linewidth]{\includegraphics[width=0.475\textwidth]{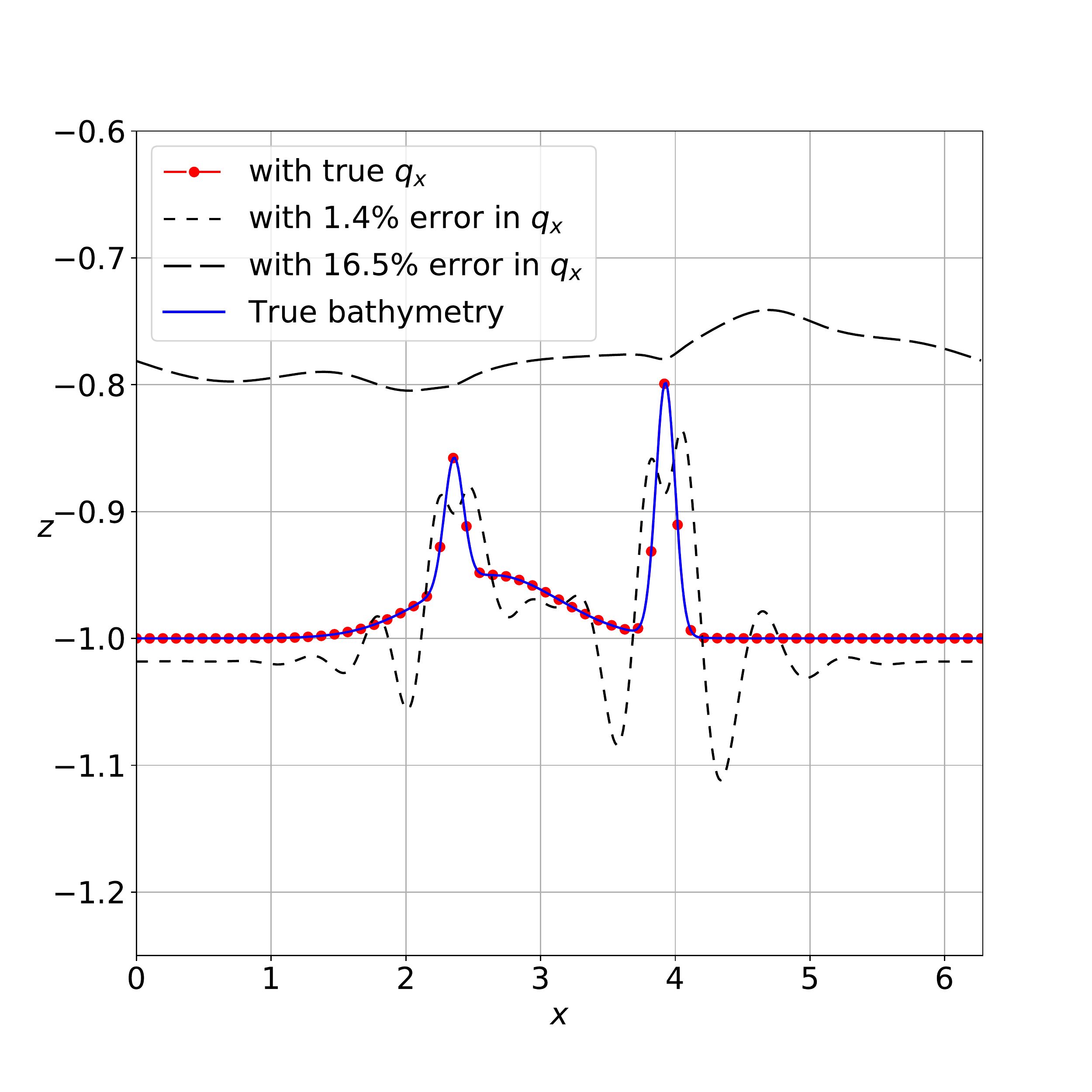}}
\caption{Bottom-boundary reconstruction using (\ref{eqn:bathymetry_min_alltime}) for regularised Boussinesq-Whitham equation (\ref{eqn:regBoussWhitham}). The black dashed lines indicate reconstruction using an erroneous value for $q_x$ with relative error percentages indicated in the legend. The construction of these inaccurate $q_x$ profiles is explained in Section \ref{sec:simultaneous}. Note even a small amount relative error can lead to inaccurate reconstruction. The error between the true profile (solid blue line) and the reconstruction using the true $q_x$ (red line with filled markers) is below machine precision. }
\label{fig:reconstruction:BoussWhitham}
\end{figure}

\begin{figure}
\centering
\subcaptionbox{Profile \ref{fig:profile_wavy} }[0.475\linewidth]{\includegraphics[width=0.475\textwidth]{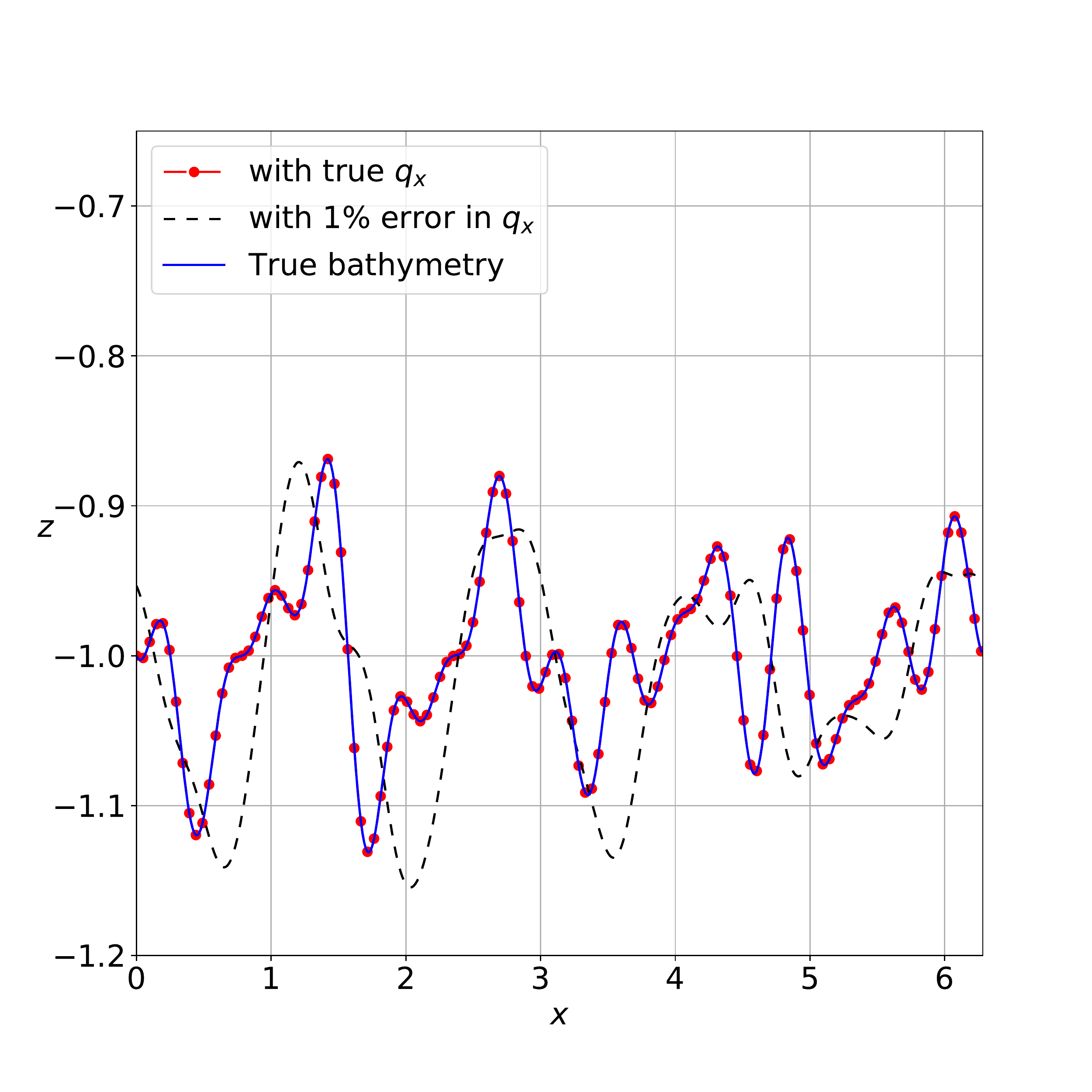}}
\subcaptionbox{Profile \ref{fig:profile_isolated}}[0.475\linewidth]{\includegraphics[width=0.475\textwidth]{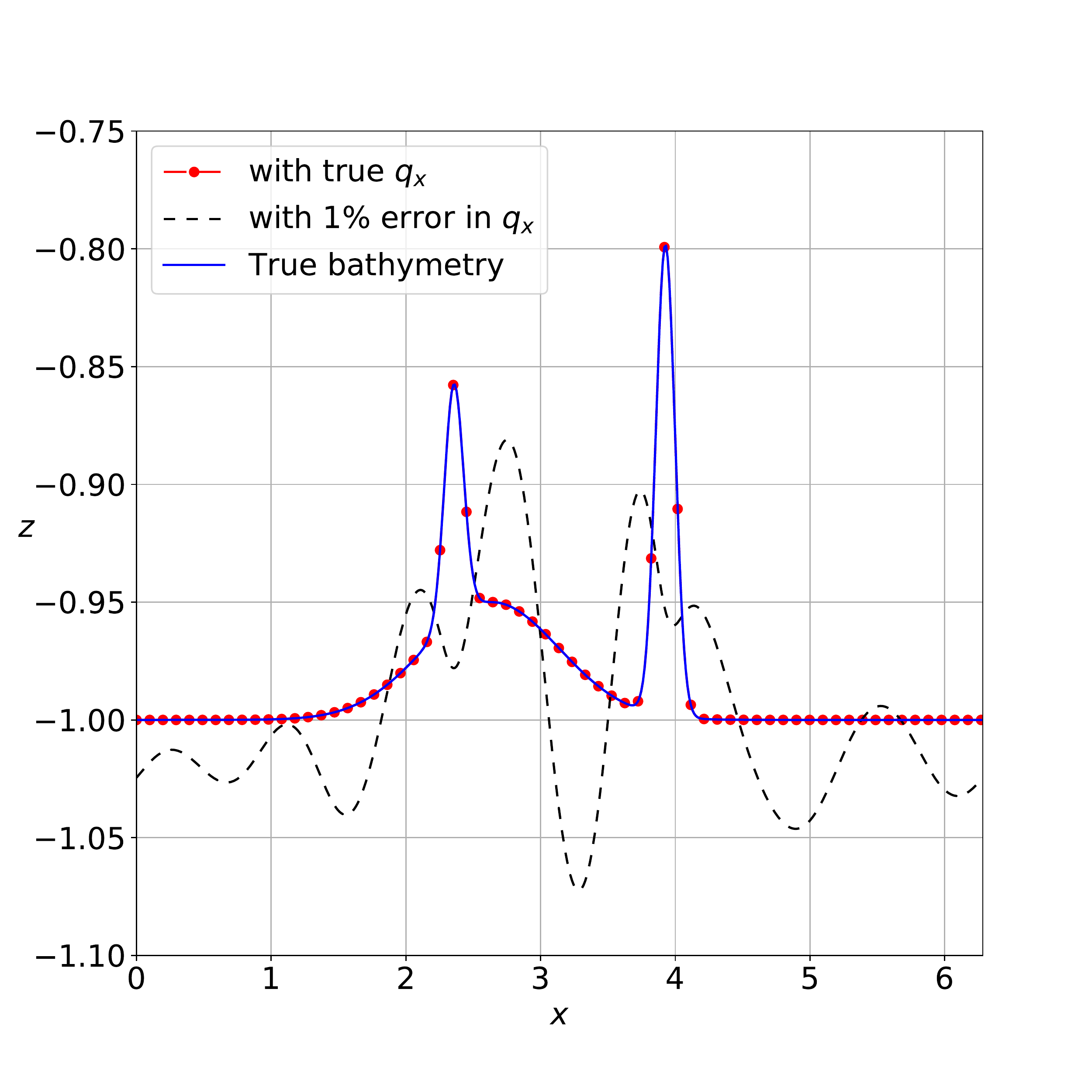}}
\caption{Bottom-boundary reconstruction using (\ref{eqn:bathymetry_min_alltime}) for regularised Boussinesq equation (\ref{eqn:regBouss}). The black dashed lines indicate reconstruction using an erroneous value for $q_x$ with relative error percentages indicated in the legend. The construction of these inaccurate $q_x$ profiles is explained in Section \ref{sec:simultaneous}. Note even a small amount relative error can lead to inaccurate reconstruction. The error between the true profile (solid blue line) and the reconstruction using the true $q_x$ (red line with filled markers) is below machine precision. }
\label{fig:reconstruction:Bouss}
\end{figure}

We now present the results of some numerical experiments investigating the inverse problem formulated in this section. First we analyse the problem of inverting the operator on the left-hand side of equation (\ref{eqn:bathymetry_EL:all_time}). To do so, we consider a particular form for the function $q_x$. Specifically we chose $q = 0.1\sin(x-t),\: x\in[0,2\pi]$ which represents a periodic travelling-wave profile for the velocity potential. We discretised this function on a uniform grid in the $x-$variable with $256$ points. For fixed $t$, using a pseudospectral method we compute the matrix representation of the operator $(\mathcal{P}q_x)\mathcal{P}^2\partial_x^2((\mathcal{P}q_x)\: \cdot)$  for both regularised Boussinesq (\ref{eqn:regBouss}) and regularised Boussinesq-Whitham (\ref{eqn:regBoussWhitham}).  This generates a $256$ by $256$ real symmetric matrix for either model. We repeated this procedure for different $t$-values, specifically for $t=2\pi j/M$ where $j=1,2,\ldots M$. This results in $M$ matrices each of size $256$ by $256$ (again for either model). $M$ matrices are summed to produce a sole real symmetric matrix of size $256\times 256$ for regularised Boussinesq and a different $256\times 256$ real symmetric matrix for regularised Boussinesq-Whitham. We have thus computed the operator on the left-hand side of (\ref{eqn:bathymetry_EL:all_time}). We then used the standard LAPACK subroutines to determine the associated eigenvalues. This procedure was repeated for different values of $M$, specifically $M=1,10,100,200$. As expected, all the eigenvalues were negative and clustered near the origin. In Figure \ref{fig:eigvals} we plot the magnitude of the sorted eigenvalues (as a function of index) for different values of $M$, indicated in the legend. We first note, that the operator when considered at a single instance of time $M=1$ is the most singular/ill-conditioned for both model equations. However the effect of taking into account additional times is indeed a form of regularisation. The eigenvalues for larger values of $M$, show an overall increase in magnitude. This is particularly so for regularised Boussinesq-Whitham which almost levels off completely. Our other model, regularised Boussinesq, shows a markedly different behaviour. Although there is some regularisation obtained for $M>1$, we still see a rapid decay in the magnitude of the eigenvalues. This is consistent with the fact that the operator for this model is in fact compact (indeed the sum of $M$ compact operators) whereas for regularised Boussinesq-Whitham the operator is only bounded. 

The rapid decay of the eigenvalues will have a strong impact on the accuracy of the bottom profile reconstruction. Note that when $\eta,\eta_t,q$ are consistent with the model equations, we have 
\[\eta_t^{(j)} = \omega^2 q^{(j)} - \mathcal P\partial_x\big( (\eta^{(j)}+\zeta)\: (\mathcal Pq_x^{(j)}) \big),\mbox{ for each }j=1,\ldots,M.\] Hence the equation we need to solve for $\zeta^*$ becomes
\[\sum_{j=1}^M(\mathcal Pq_x^{(j)}) \: \mathcal P^2\partial_x^2 \left( (\mathcal Pq_x^{(j)}) \: \zeta^*\right) = \sum_{j=1}^M (\mathcal Pq_x^{(j)}) \: \mathcal P^2\partial_x^2 \left( (\mathcal Pq_x^{(j)}) \: \zeta \right),\] 
which provides the true solution \emph{if} we can invert the matrix on the left side of the equality. If the data $\eta,\eta_t,q$ provided are not consistent or only approximately satisfy the necessary relationship between themselves, then these small errors will inevitably get magnified when reconstructing the bottom topography. Hence the decay of the eigenvalues determines the relative accuracy to which we require the provided data to be consistent with the underlying model.

{\remark The above analysis was for a fixed spatial resolution, \emph{i.e.} fixed number of grid points in $x$. If the resolution in $x$ were increased,  one may require higher values of $M$ to achieve a similar amount of regularisation. For regularised Boussinesq, fixing $M$ and increasing the resolution in $x$, will still result in small eigenvalues and cause the operator on the left-hand side of (\ref{eqn:bathymetry_min_alltime}) to be ill-conditioned.  Thus we are assuming the bottom-profile to be recovered is well approximated with the chosen number of grid points. Equivalently we assume the bottom-profile is not too rough. This is not unexpected in these kinds of ill-conditioned inverse problems. However, Figure \ref{fig:eigvals} indicates the regularised Boussinesq-Whitham model behaves quite differently under this regularisation scheme. Recall the associated operator here is not compact. The primary source of ill-posedness comes from the zeros of $\mathcal Pq_x$. When the zeros of $\mathcal Pq_x$ are not stationary in time, taking larger $M$ seems to help for this model. We observed similar behaviour at higher resolutions in $x$.\\}

We continue our experiments with reconstruction by considering examples using data taken from forward simulations of the model equations themselves. In Figures \ref{fig:reconstruction:BoussWhitham} and \ref{fig:reconstruction:Bouss} we consider the recovery of the bottom topography profile $\zeta$ given $\eta,\eta_t,q$ from a forward simulation of the respective model equations. For each model equation, we consider the recovery of a bottom profile which consists of a number of a sine waves called Profile 1 (Figure \ref{fig:profile_wavy}) or a localised bottom feature called Profile 2 (Figure \ref{fig:profile_isolated}). In all cases, a solid blue line indicates the true bottom topography and the red line with filled-circle markers indicates our recovered profile when using data obtained at $200$ different times from the forward simulation. For a given time instance $t_j$, for which we seek to compute the left-hand side operator in equation (\ref{eqn:bathymetry_EL:all_time}), we recorded the solution $(\eta,q)$ to the forward problem from a numerical simulation of the model equations. This was done in such a way that the time-derivative $\eta_t$ could be computed using a fourth-order accurate five-point finite-difference stencil in $t$. This stencil lead to a time-derivative of sufficient accuracy so as to enable recovery and was, up to machine precision, equivalent to computing the right-hand side of equation (\ref{eqn:omegaPa}). We define two relative errors for the recovered bottom-boundary as 
\begin{align}\label{errors}
	\mbox{E}_b = \frac{\|(1+\zeta_r)-(1+\zeta)\|_2}{\|1+\zeta\|_2},\quad \mbox{E}_p =  \frac{\|\zeta_r-\zeta\|_2}{\|\zeta\|_2}
\end{align}
where $\zeta_r$ is the reconstructed profile, $\zeta$ is the true profile, $\mbox{E}_b$ denotes the relative error in the actual depth (including the baseline $z=-1$) and $\mbox{E}_p$ is the relative error for the deviation from $z=-1$. The value of $\mbox{E}_b$ for the recovered bottom-boundary with data $(\eta,\eta_t,q_x)$ consistent with the forward model (for either regularised Boussinesq or regularised Boussinesq-Whitham) is approximately $10^{-12}$, indicating recovery is possible. We also used a second-order accurate three-point stencil to compute $\eta_t$. Despite a higher $\mbox{E}_b$ (around $10^{-8}$) we could still reconstruct the bottom-boundary.

On the other hand, using a value of $q_x$ with a small amount of error, gave rise to considerable error in the recovered bottom topographies. The reconstructions using erroneous values of $q_x$ are shown in Figures \ref{fig:reconstruction:BoussWhitham} and \ref{fig:reconstruction:Bouss} by long-dashed black lines. The recovered profiles with even a $1\%$ $L^2-$norm relative-error in $q_x$ (but no error in $\eta,\eta_t$) resulted in  $\mbox{E}_b$ between $4-7\%$. Admittedly recovering bottom topographies through this approach is prone to noise and error, despite our regularisation. For the case of regularised Boussinesq, we do not consider reconstruction corresponding to $q_x$ values with error larger than a few percent since it gave such poor results that one could not properly visualise all recovered profiles on a single plot. In this sense, regularised Boussinesq-Whitham is arguably better; a reflection of the slower decay in the magnitude of the eigenvalues for this model. We did try other Tikhonov-based regularisation schemes. However these methods needed considerable tuning of the regularisation parameter which was particular to each bottom-profile and we were unable to find a systematic way to do so.

\section{The observer model for velocimetry}\label{sec:observers}
In this section we consider the problem of determining the velocity $q_x$ of the fluid from measurements of the surface elevation $\eta$ when the bottom boundary $\zeta$ is known and $\|\zeta\|_\infty$ is finite. To achieve this we employ the observer framework \cite{luenberger1971introduction}. An observer system is a dynamical system  which is constructed in relation to another dynamical system when given partial knowledge of the state of the latter dynamical system. Let $y$ represent the state of a dynamical system that evolves according to $y_t = f(y)$. An observer for this system is $\tilde y_t = \tilde f(\tilde y,Oy)$ where $\tilde f$ is a modification of $f$ and $O$ is an operator with a null-space representing the fact that we only have knowledge of a part of the state $y$. The  goal then is to design $\tilde f$ such that $\tilde y \to y$ as $t\to\infty$. Hence by simulating the observer model, we may retrieve knowledge of the full state. 

Let $q,\eta$ represent the solution to (\ref{eqn:omegaPa}-\ref{eqn:omegaPb}) for some unknown initial condition. These represent the true state. The observer variables are denoted by $\tilde q,\tilde \eta$ and we propose they satisfy 
\begin{align}\label{eqn:observer_NLa}
	\tilde \eta_t &= \omega^2 \tilde q - \mathcal P\partial_x \big( (\tilde \eta + \zeta)\:(\mathcal P\tilde q_x)\big)-\lambda(\tilde \eta - \eta),\\
	\tilde q_t &= -\tilde \eta - \frac{1}{2}\left(\mathcal P\tilde q_x \right)^2 - \nu(\tilde \eta - \eta).\label{eqn:observer_NLb}
\end{align}
Here $\lambda,\nu$ are the observer parameters. In this section we assume we have access to the function $\eta$, or equivalently we are provided the surface deviation for all $x,t$ from measurements. The problem at hand is then to recover the velocity $q_x$ (which cannot be measured) by solving (\ref{eqn:observer_NLa}-\ref{eqn:observer_NLb}) with arbitrary initial conditions for $\tilde\eta,\tilde q$, assuming we know the bottom-profile $\zeta$. Thus we seek a rationale for choosing $\lambda,\nu$ to ensure the required convergence, namely $\tilde q_x\to q_x$ as $t\to\infty$.  At the end of this section, we present numerical simulations that validate this choice. Our numerical simulations employ the surface deviation $\eta(x,t_j)$ for a discrete set of times $t_j$.

Let $q^e = \tilde q - q$, $\eta^e = \tilde \eta - \eta$, represent the errors in the potential and free-surface respectively. We first state a theorem for the error associated with the linear constant-coefficient part of equations (\ref{eqn:observer_NLa}-\ref{eqn:observer_NLb}). 
{\theorem{\label{thm:cc:decay}
	For any positive number $d$, there exist real numbers $\lambda$ and $\nu$ such that the solution to 
	\begin{align}\label{eqn:obsever:cca}
		\eta^e_t &= -\lambda \eta^e + \omega^2  q^e,\\
		q^e_t &= -(1+\nu)\eta^e,\label{eqn:obsever:ccb}
	\end{align}
	with periodic boundary conditions satisfies \[\|\eta^e\|_2 \leq C e^{-dt},\quad \|q_x^e\|_2 \leq C e^{-dt},\]
	where $C$ is a constant that depends on the initial condition.

	\begin{proof}
		The proof is straightforward and follows from the Fourier series solution to the equations. Indeed the equations for $\eta^e,q^e$ are equivalent to \[\eta_{tt}^e + \lambda \eta_t^e + (1+\nu)\omega^2 \eta^e = 0,\]
		which has solutions that look like $\exp(ikx+p_k t)$ where \[p_k^2 + \lambda p_k + (1+\nu)\omega^2(k) = 0\Rightarrow p_k = -\frac{\lambda}{2} \pm \sqrt{\frac{\lambda^2}{4}-(1+\nu)\omega^2(k)}\]
		Note for both choices of regularised Boussinesq (\ref{eqn:regBouss}) and regularised Boussinesq-Whitham (\ref{eqn:regBoussWhitham}), $\omega^2(k)\geq 0$, vanishes only when $k=0$ and is an increasing function of $|k|$. Hence setting $\lambda = 2d$, we obtain the required decay rate by ensuring the term under the square root is negative. This can be done by choosing $\nu$ so that
		\begin{align}
			1+\nu >\left\{ \begin{array}{ll} d^2\:\dfrac{1+\mu^2/2}{1+\mu^2/6},\quad &\mbox{for (\ref{eqn:regBouss})},  \\ \\ d^2\:\dfrac{\mu}{\tanh(\mu)},\quad &\mbox{for (\ref{eqn:regBoussWhitham})}.  \end{array} \right.
		\end{align}
		The statement of the theorem then follows by writing the Fourier series solution and computing an estimate for sufficiently smooth initial data.
	\end{proof}
}}
It follows from the above theorem that we can design the observer to guarantee exponentially fast error decay at any desirable rate for the constant coefficient part of (\ref{eqn:observer_NLa}-\ref{eqn:observer_NLb}). Note that the error in the mean mode of $q$ does not decay. We observe this from the dispersion relation in the proof of the theorem where $p_k=0,-\lambda$ when $k=0$. We also note more directly from the integral of linearised version of equations (\ref{eqn:observer_NLa}-\ref{eqn:observer_NLb})
\[\partial_t \int \eta^e = -\lambda \int \eta^e,\quad\quad  \partial_t\int q^e = -(1+\nu)\int \eta^e,\] and note that the mean mode of $\eta^e$ vanishes exponentially, but the initial error in $q^e$ is never eliminated. For the linear observer equations, one can only recover the velocity and not the potential. This is physically reasonable.

\subsection{Decay of error in the linear equations}\label{sec:decay_linear}
The previous theorem not only gave us the required decay rate but it also provided (in principle) a solution expression for the constant coefficient part of the linear operator on the right-hand side of  (\ref{eqn:observer_NLa}-\ref{eqn:observer_NLb}) and a representation for the associated semi-group $e^{\mathcal L t}$. On the other hand, the linear operator appearing in (\ref{eqn:observer_NLa}-\ref{eqn:observer_NLb}) is in fact a variable-coefficient operator due to the presence of $\zeta$ and we would like to estimate the decay rate of solutions associated to the entire linear part. We claim that the linear equation given by
\begin{align}\label{eqn:observer_linear:zetaa}
	 \eta_t^e &= \omega^2  q^e - \mathcal P\partial_x \left( (\mathcal P q_x^e)\: \zeta\right)-\lambda \eta^e,\\
	 q_t^e &= -(1+\nu) \eta^e,\label{eqn:observer_linear:zetab}
\end{align}
possesses a solution which can be given in terms of a continuous one-parameter family of solution operators. This follows from the fact that 
\[- \mathcal P\partial_x \left( (\mathcal P\tilde q_x)\: \zeta\right),\] is a bounded self-adjoint operator acting on $\tilde q$, for both regularised Boussinesq and regularised Boussinesq-Whitham, and employing standard results in the perturbation theory of operators \cite{kato2013perturbation}. Next we consider (\ref{eqn:observer_linear:zetaa}-\ref{eqn:observer_linear:zetab}) in integral form as follows
\begin{align}
	\begin{pmatrix}
	\eta^e \\ q^e
	\end{pmatrix}
	&=
	e^{\mathcal L t}
	\begin{pmatrix}
	\eta^e(x,0) \\ q^e(x,0)
	\end{pmatrix}
	+ 
	\int_0^t e^{\mathcal L (t-s)}
	\begin{pmatrix}
		-\mathcal P \partial_x\left(\zeta  \: (\mathcal P q_x^e) \right)
		\\
		0
	\end{pmatrix}ds.
\end{align}
If we limit ourselves to initial data with zero average, then Theorem \ref{thm:cc:decay} affords a uniform decay-rate and thus we have
\begin{align}
	\|y\| \leq e^{-\frac{\lambda t}{2} } \|y_0\| + \int_0^t e^{-\frac{\lambda}{2}(t-s) }\frac{\|\zeta\|_\infty}{\mu^j}\|y\|ds,
\end{align}
where $j=2$ for (\ref{eqn:regBouss}) and $j=1$ for (\ref{eqn:regBoussWhitham}). Here $\|y\| = \sqrt{\|\eta^e\|_2^2 + \|q_x^e\|_2^2}$ is the norm of the solution for equations (\ref{eqn:observer_linear:zetaa}-\ref{eqn:observer_linear:zetab}). The above claim follows from the triangle inequality and the fact that the induced $L^2-$norm of $\mathcal P$ is bounded by $1$ and the norm of $\mathcal P\partial_x$ depends on the particular choice of shallow-water model.  We now employ a simple Gr\"onwall argument to conclude that 
\[\|y\| \leq e^{-m t}\|y_0\|,\quad m = \frac{\lambda}{2} - \frac{\|\zeta\|_\infty}{\mu^j}.\]
The zero average condition for $q^e$ and $\eta^e$ is clearly maintained by (\ref{eqn:observer_linear:zetaa}-\ref{eqn:observer_linear:zetab}). The Gr\"onwall's inequality we employed above follows from the following theorem due to Zadiraka \cite{zadiraka1968issledoavanie} upon assuming $b=0$.
{\theorem{\label{thm:gronwalls}
Let $u(t)$	 be a continuous function that satisfies 
\[|u(t)| \leq |u(0)|e^{-\alpha t} + \int_0^t e^{-\alpha(t-s)}( c|u(t)| + b)ds,\]
where $a,b,c$ are all positive constants then 
\[|u(t)| \leq |u(0)|e^{-(\alpha-c)t} + \frac{b}{\alpha - c} (1 + e^{-(\alpha-c)t}).\]
}}
{\remark{The quoted result due to Zadiraka appears in the literature with the term $e^{-\alpha t}|u(0)|$ in the result, as opposed to $e^{-(\alpha-c)t}|u(0)|$. However we could neither find the original reference nor prove the result ourselves. Hence we choose to refer to the result above, which we could prove. }}

We conclude then that the linear equations given in (\ref{eqn:observer_linear:zetaa}-\ref{eqn:observer_linear:zetab}), like the constant coefficient counterpart, also satisfy an exponential decay bound when $\lambda$ is taken sufficiently large. The error in the non-zero modes, for the linear observer problem, decays monotonically. This follows from the error bounds derived above though one may also reach the same conclusion via energy considerations for the linear observer problem, see Appendix \ref{app:energy}.

Note that the nonlinear observer equations (\ref{eqn:observer_NLa}-\ref{eqn:observer_NLb}) do not preserve the zero-mode for the potential $q$, neither for regularised Boussinesq nor for regularised Boussinesq-Whitham. Thus we do not expect $\tilde q\to q$ as $t\to\infty$. The nonlinear terms in the observer model (\ref{eqn:observer_NLa}-\ref{eqn:observer_NLb}) and shallow-water model (\ref{eqn:omegaPa}-\ref{eqn:omegaPb}) are the same and are locally Lipschitz functions of their arguments. One may then anticipate a reduction in the error $q_x^e,\eta^e$ for finite time and $\lambda$ large enough. However we do not prove the long-time convergence for the nonlinear problem which would require a long-time wellposedness result for the underlying models. Nonetheless, our numerical simulations of (\ref{eqn:observer_NLa}-\ref{eqn:observer_NLb}) indicate the error in velocity $ q_x^e$ does indeed vanish asymptotically in time.

\subsection{Choosing the observer parameters in practice}\label{sec:observer:practice}
We now present our scheme to set the observer parameters $\lambda$ and $\nu$ for general bottom boundary profiles. Consider the linear equation for the error terms $q^e$ and $\eta^e$
\begin{align}
	 \eta_t^e &= \omega^2  q^e - \mathcal P\partial_x \left( (\mathcal P q_x^e)\: \zeta_c\right)-\lambda \eta^e,\\
	 q_t^e &= -(1+\nu) \eta^e,
\end{align}
where $\zeta_c$ is a constant. Rewriting this equation in terms of $\eta^e$ alone we have
\begin{align}
	\eta_{tt}^e + \lambda \eta_t^e + (1+\nu)(\omega^2 \eta^e -\mathcal P \partial_x ( (\mathcal P \eta_x^e)\: \zeta_c )) = 0,
\end{align}
which has solutions of the form $e^{ikx+p t}$ where $p$ satisfies
\[p^2 + \lambda p + (1+\nu)(\omega(k)^2 + k^2\mathcal P(k)^2\zeta_c) = 0\] 
\[\Rightarrow p = -\frac{\lambda}{2} \pm \sqrt{\frac{\lambda^2}{4} - k^2(1+\nu)\left(\frac{\omega(k)^2}{k^2} + \zeta_c \mathcal{P}(k)^2\right)}\]
As in the case of Theorem \ref{thm:cc:decay} we need to ensure the term under the square root is negative for all $k$. The coefficient of $(1+\nu)$ is an increasing function of $k$ for both choices of $\omega^2,\mathcal P$ in (\ref{eqn:regBouss}) and (\ref{eqn:regBoussWhitham}). Hence it suffices to ensure the appropriate sign under the root for $k=1$. If $\zeta_c$ is positive then the parameter choice of Theorem \ref{thm:cc:decay} is sufficient here too. This suggests that for a general profile $\zeta$, a suitable choice for $\zeta_c = \min[0,\zeta]$. In keeping with the no-island condition, we assume $\zeta>-1$ then the energy corresponding to solutions of equations (\ref{eqn:observer_linear:zetaa}-\ref{eqn:observer_linear:zetab}) is positive and monotonically decreases (see Appendix \ref{app:energy}). Hence we have for the constant coefficient equations a
\begin{align}
\mbox{decay rate of }\lambda/2\mbox{ when } (1+\nu)> \left.\frac{\lambda^2}{4(\omega(k)^2 + \zeta_c\mathcal P(k)^2)}\right|_{k=1}\: \zeta_c=\min[0,\zeta].\label{decay_rate}
\end{align}
As shown in our numerical experiments, this choice seems to give the decay rate of $\lambda/2$ for the variable coefficient equation (\ref{eqn:observer_linear:zetaa}-\ref{eqn:observer_linear:zetab}) as well as the fully nonlinear observer model (\ref{eqn:observer_NLa}-\ref{eqn:observer_NLb}). 

\begin{figure}
\centering
\subcaptionbox{Error decay for profile \ref{fig:profile_wavy}}[0.475\linewidth]{\includegraphics[width=0.475\textwidth]{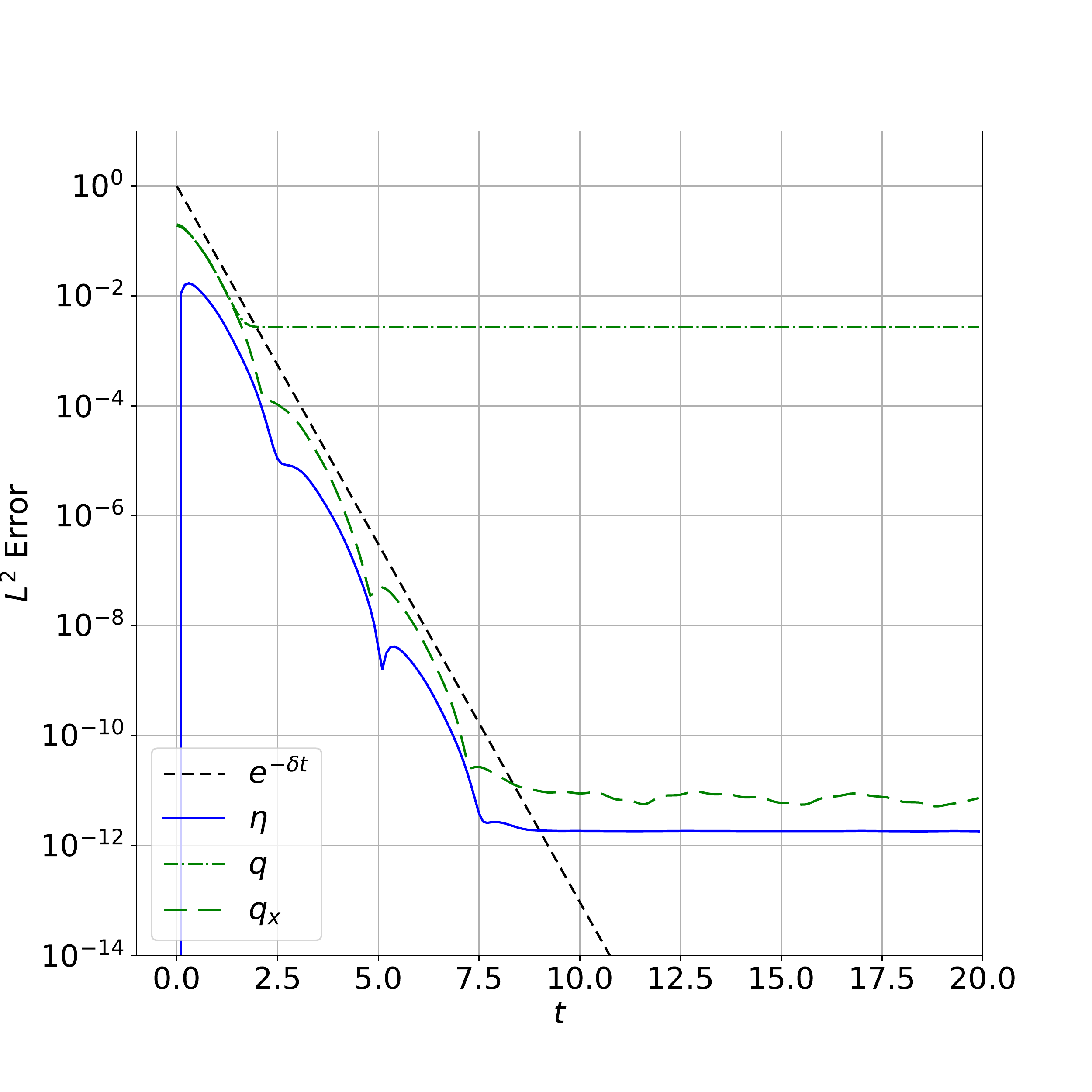}}
\subcaptionbox{Error decay for profile \ref{fig:profile_isolated}}[0.475\linewidth]{\includegraphics[width=0.475\textwidth]{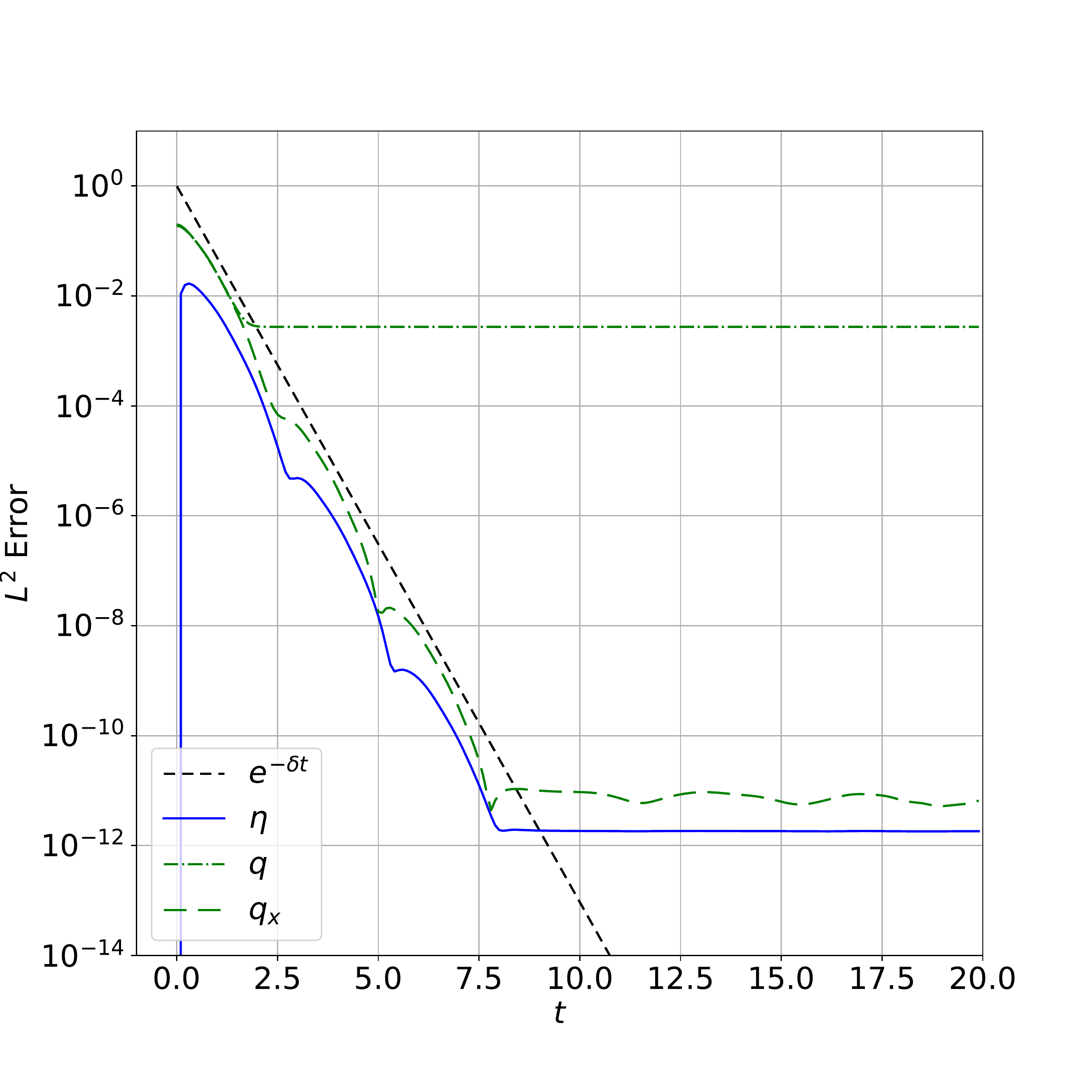}}
\caption{Decay of error in solution for the observer model corresponding to regularised Boussinesq-Whitham (\ref{eqn:regBoussWhitham}) with observer parameters $\lambda=6$ and $\nu=14$. This results in a linear decay rate of $\lambda/2=3$ which is indicated in the short-dashed black line. The error for $\eta$ and $q_x$ in the full nonlinear observer problem follow almost exactly the predicted linear decay rate. Shown in the green dashed-dot line is the error in $q$ which saturates to a non-zero value.}\label{fig:observer:BoussWhitham}
\end{figure}

\subsection{Numerical experiments}\label{sec:experiments:observer}

\begin{figure}
\centering
\subcaptionbox{Error decay for profile \ref{fig:profile_wavy}}[0.475\linewidth]{\includegraphics[width=0.475\textwidth]{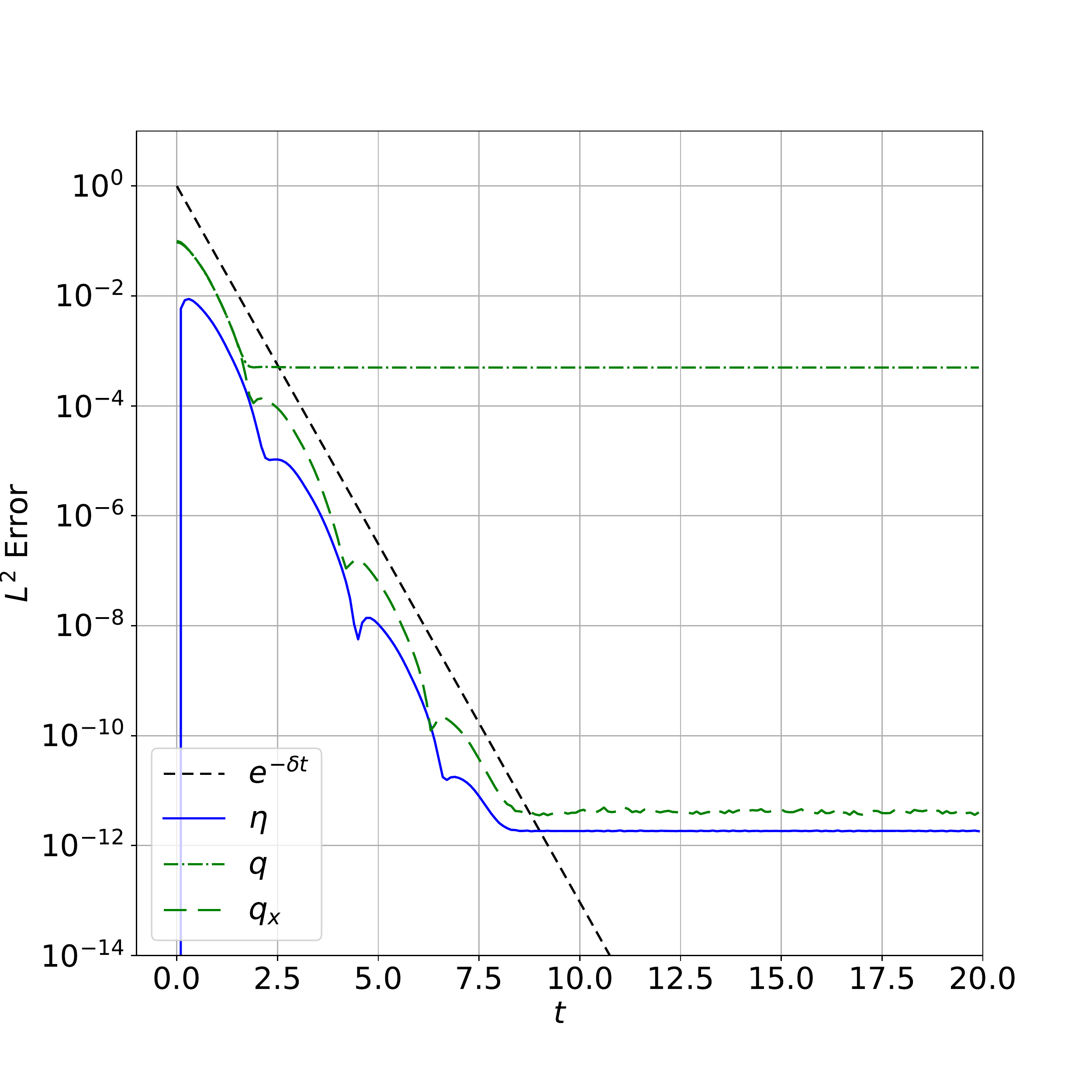}}
\subcaptionbox{Error decay for profile \ref{fig:profile_isolated}}[0.475\linewidth]{\includegraphics[width=0.475\textwidth]{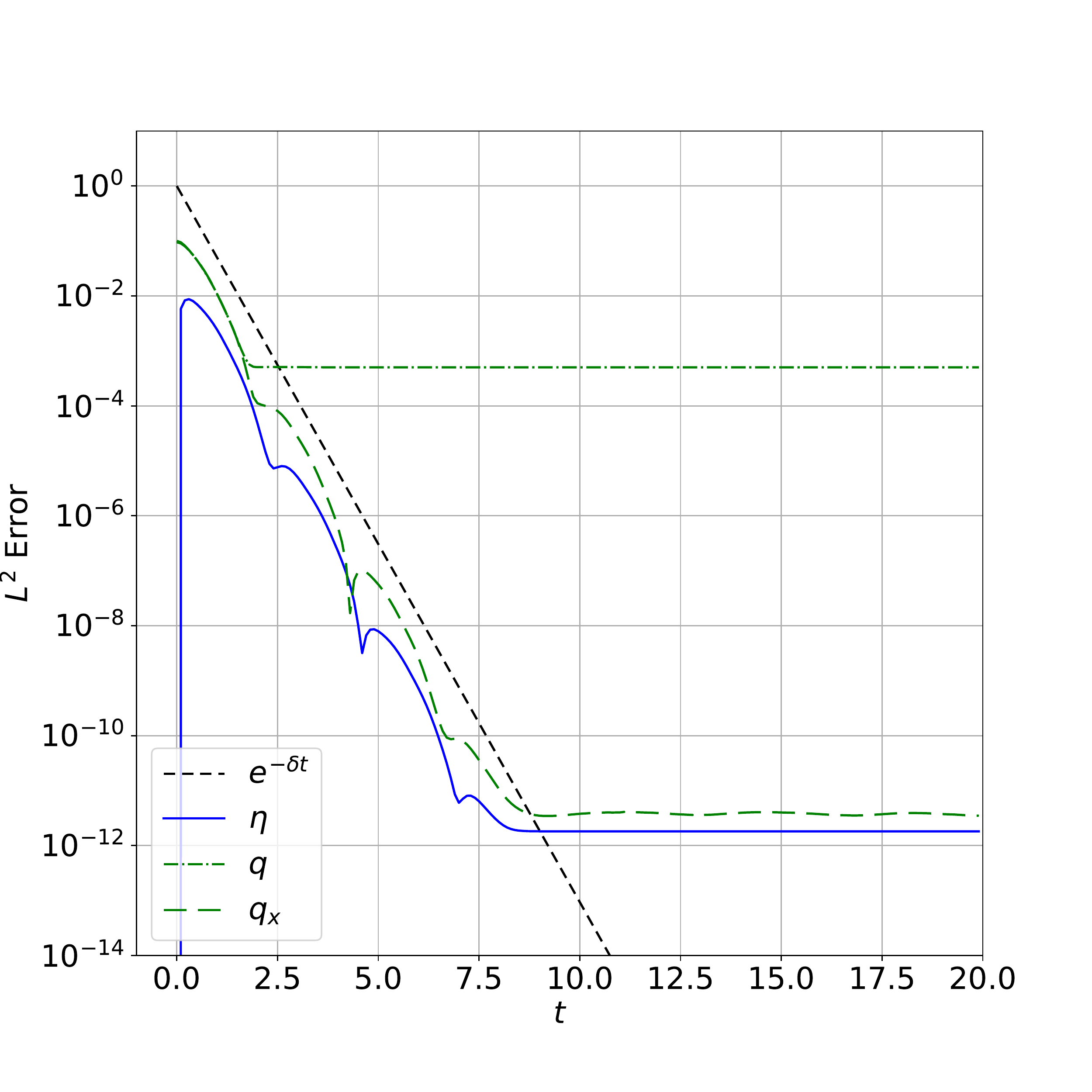}}
\caption{Decay of error in solution for the observer model corresponding to regularised Boussinesq (\ref{eqn:regBouss}) with observer parameters $\lambda=6$ and $\nu=14$. This results in a linear decay rate of $\lambda/2=3$ which is indicated in the short-dashed black line. The error for $\eta$ and $q_x$ in the full nonlinear observer problem follow almost exactly the predicted linear decay rate. Shown in the green dashed-dot line is the error in $q$ which saturates to a non-zero value.}\label{fig:observer:Bouss}
\end{figure}

To evolve the solution of the observer equations (\ref{eqn:observer_NLa}-\ref{eqn:observer_NLb}) forward in time, we require the solution of the model shallow-water equations (\ref{eqn:omegaPa}-\ref{eqn:omegaPb}), or more specifically $\eta(x,t)$. Consequently, we solve both equations simultaneously with the observer model coupled to the original equations. Since we employ periodic boundary conditions, we represent all functions $\eta,q,\tilde \eta,\tilde q$ as well as the bottom-boundary using Fourier series and employ a pseudospectral method, with $2/3$-method for de-aliasing, to time-evolve the full system of (four) equations. In all the simulations, in this section and the next, we used $512$ grid points in the $x-$direction. We used the standard explicit fourth-order Runge-Kutta scheme with a time-step $\Delta t = 10^{-3}$ non-dimensional time units for the model shallow-water equations (\ref{eqn:omegaPa}-\ref{eqn:omegaPb}) and recorded $\eta$ at each time-step. We also used fourth-order RK4 for the observer model but with a time-step of $2\times 10^{-3}$ to account for the fact that the solution $\eta$ is supplied to the observer model at a fixed rate and since the RK4 scheme requires the value of the vector field at intermediate time steps. Additionally we conducted experiments where the solution to the model equations was recorded every $m$ time steps and hence the observer model was time evolved with a time step of $2m\Delta t$. In all cases we obtained the same qualitative behaviour. In the present manuscript we only show results for $m=1$ and $\Delta t= 10^{-3}$.

Recall we defined the errors in the velocity potential and surface deviation as $q^e = \tilde q - q$ and $\eta^e = \tilde \eta - \eta$, where the tilde variables correspond to the solution of the observer model and $\eta,q$ are the solution of our model shallow-water equations. Figures \ref{fig:observer:BoussWhitham} and \ref{fig:observer:Bouss} show the $L^2$ norms of $q^e$ and $\eta^e$ as a function of time for the two different models, regularised Boussinesq-Whitham and regularised Boussinesq respectively. In each plot, we also show the error in the velocity $q_x^e = \tilde q_x - q_x$ as well as the predicted decay rate according to our observer design, which in this case corresponds to $\lambda/2$. For all cases, we took $\lambda=6$ and $\nu = 14$. This choice is consistent with the design suggested in (\ref{decay_rate}) for either of the bottom profiles we considered. We note that the errors for the full nonlinear observer problem follow the predicted decay rate very closely. We conjecture that $\lambda$ was sufficiently large to overcome any possible growth in the error due to the nonlinearity. 
We also note that the error in the velocity potential does not decay to zero, as expected. The decay in the error did not depend substantially on which of the the bottom profiles we considered. Figures \ref{fig:observer:BoussWhitham} and \ref{fig:observer:Bouss} indicate the errors in the surface deviation and velocity do not decrease to machine precision. The value to which they saturate depends on the frequency with which the data $\eta(x,t)$ of the original model is recorded. As we decreased $\Delta t$, the final asymptotic-in-time errors decreased. We conclude this section by emphasising the main upshot of our numerical simulations: when the bottom-boundary $\zeta$ is known, it is possible to recover the velocity of the fluid in these shallow-water models given only measurements of the surface deviation. The error in the recovered velocity depends on the frequency of the measurements.

\section{Simultaneous velocity-recovery and bathymetry}\label{sec:simultaneous}
Previously we considered the problem of recovering $\zeta$ given $\eta,\eta_t,q_x$ (Section \ref{sec:bathymetry_reconstruct}) and then the problem of recovering $q_x$ given $\zeta$ and $\eta$ (Section \ref{sec:observers}). In the present section we combine these two inverse problems into one, and propose a method to recover $\zeta$ from measurements of $\eta$ alone. This will entail recovering $q_x$ as well. Henceforth we assume we know $\eta$ (as a function of $x$ and $t$). The recovery of velocities and the bottom-boundary involves carefully selecting the observer parameters $\lambda,\nu$. To do so, we need to analyse the errors associated with the methods described in Sections \ref{sec:bathymetry_reconstruct} and \ref{sec:observers}.

This section consists of three subsections. In the first subsection, we describe the overall reconstruction algorithm. Subsequently, we obtain error estimates that motivate our particular choice for the parameters $\lambda,\nu$. As the reconstruction algorithm consists of two stages, we first derive error estimates for both stages independently and then combine them. The upshot is a condition on the observer parameters. Finally we conclude with some example reconstructions, \emph{i.e.} starting with an initial guess for the bottom-profile and measurements of the free surface $\eta$ as functions of $x$ given at specific instants of time, we improve our initial estimate to determine the true bottom-profile.

\subsection{The reconstruction algorithm}\label{sec:reconstruction_algo}
 The basic strategy is as follows. Suppose we have an \textit{a priori} estimate for $\zeta^e$, or in other words a reasonable but not necessarily accurate guess for the bottom boundary. We employ this initial guess in the observer equations (\ref{eqn:observer_NLa}-\ref{eqn:observer_NLb}) to arrive at an estimated $\tilde q_x$ which we record at multiple instances of time. We then employ the estimated velocity $\tilde q_x$ in the reconstruction equation (\ref{eqn:bathymetry_EL:all_time}) to update our estimate of the bottom boundary profile. 

We now present the main contribution of this manuscript. The following is our proposed algorithm for bathymetry using surface-wave measurements alone.
\begin{enumerate}
	\item Suppose we are given $\eta(x,t)$, the true surface deviation, and a reasonable guess for the bottom boundary $\tilde \zeta$. We assume the constant function $\tilde \zeta(x)=\zeta_c$.
	\item Pick a small number $\epsilon$ such that $\epsilon\ll \mu^2$. Choose $\lambda$ and $\nu$ such that
	\[\frac{1+\nu}{\lambda} = \epsilon,\quad \mbox{ and }\quad  (1+\nu)> \left.\frac{\lambda^2}{4(\omega(k)^2/k^2 + \zeta_c\mathcal P(k)^2)}\right|_{k=1}. \]
	Note the second condition above ensures the error in the observer model decays. The motivation for the first condition is given in the next subsection.
	\item Run the observer model (\ref{eqn:observer_NLa}-\ref{eqn:observer_NLb}) with these parameters (using any reasonable initial condition; we employ trivial initial conditions) for a time till the predicted error $e^{-\frac{\lambda t}{2} }$ is sufficiently small and record $\tilde q_x$ from this simulation for a large number of time instances (say $200$).
	\item Using $\eta,\eta_t$ and $\tilde q_x$ at the different instances of time, solve the reconstruction equation (\ref{eqn:bathymetry_EL:all_time}) to recover the bottom boundary.
\end{enumerate}

Note our proposed method involves two steps: an observer problem followed by a reconstruction step. This implies the estimated velocity provided by the observer problem must be sufficiently accurate so that the resultant $\zeta$ we reconstruct is accurate. We emphasise that we do not iterate these two steps. The ill-conditioned nature of the matrix inversion involved in solving equation (\ref{eqn:bathymetry_EL:all_time}) caused small errors to amplify when we implemented an iterative scheme. This forced errors in both the velocity estimate from the observer problem and reconstructed bottom boundary profile to rapidly grow.

The requirements on $\lambda$ and $\nu$ imply that $\lambda$ cannot be too large. Indeed $\lambda$ is typically less than or equal to $\epsilon$. Hence, despite the fact that we were free to choose almost any $\lambda$ in the observer problem, when attempting to reconstruct the bottom-profile, we cannot impose too large a decay rate on the error; the surface deviation must be assimilated slowly into the observer equations so that the resulting velocity is accurate. We recall the examples discussed in Section \ref{sec:bathy_reconstruct:examples} which involved reconstruction from erroneous $q_x$ values. Such velocities were obtained from an observer problem (with an estimated $\tilde\zeta$) using `large' values of $\lambda$, namely those which did not satisfy the requirements detailed above. In such a case, the error in the estimated velocity was too large and this lead to the poor reconstruction of the bottom boundary profile.

{\remark Another noteworthy point is that $\nu$ is negative when $1+\nu = \epsilon \lambda$. Indeed it is close to $-1$. When $\nu$ is precisely $-1$ then the observer equation for the velocity potential is decoupled from the $\eta$ equation
\[\tilde q_t = -\eta - \frac{1}{2}\left(\mathcal P \tilde q_x\right)^2.\]
Thus the velocity potential is purely driven by the measured $\eta$. It is not evident that the error due to the mismatch in the initial condition reduces over time. Indeed it is not evident that there is any decay in $q_x^e$ whatsoever. Moreover the equation is similar to a forced inviscid Burgers equation and we suspect it is liable to produce high derivatives. The case for $\nu\to -1^+$ leads to eventual decay in the error but on a very long time-scale given by $2/\lambda$ as well as some amount of dispersive smoothing for high wave numbers as evident from the dispersion relation (see Theorem \ref{thm:cc:decay}).}

\subsection{Error estimates} We now present arguments that motivate the condition that $(1+\nu)/\lambda$ be taken as small as possible. The argument follows from error estimates, for both regularised Boussinesq and regularised Boussinesq-Whitham, in either of the two stages of our reconstruction algorithm, namely (i) reconstructing the bottom-profile $\zeta$ given $\eta(x,t_j),\eta_t(x,t_j)$ and $q_x(x,t_j)$ for $j=1,2,\ldots,M$ (as detailed in Section \ref{sec:bathymetry_reconstruct}), (ii) estimating the velocity $q_x$ given $\eta(x,t)$ and an initial guess for the bottom-profile $\tilde\zeta$ using the observer framework (as detailed in Section \ref{sec:observers}). Finally we combine the estimates to determine how the initial error in the bottom-profile may be reduced.

\subsubsection{Error estimates for bottom-profile reconstruction}
Consider equation (\ref{eqn:bathymetry_EL}) written here in terms of the errors 
$q^e=\tilde q - q$, $\eta^e=\tilde \eta - \eta$ and $\zeta^e =  \zeta^* - \zeta$
\begin{align}\label{eqn:error_zeta}
	(\mathcal P\tilde q_x) \: \mathcal P^2\partial_x^2 \big( (\mathcal P\tilde q_x) \: \zeta^e\big) = (\mathcal P\tilde q_x)\: \mathcal P\partial_x \big(\omega^2  q^e - \mathcal P\partial_x\left( (\mathcal P  q_x^e ) \: (\eta + \zeta) \right)  \big),
\end{align}
where the tilde variables denote our current estimate of the observer model and $q,\eta,\zeta$ represent the true values. This equation indicates the error in our estimation of $\zeta$ is proportional to the error in $q$. Note the operator that appears on the right-hand side, acting on $q^e$, is in fact the right-hand side of equation (\ref{eqn:omegaPa}): the DNO for the shallow-water model with the \emph{true} values of $\zeta$ and $\eta$. However, as seen in the previous section, we only have decay in $q_x^e$ and not $q^e$. Hence we consider the right-hand side as an operator on $q_x^e$.

For regularised Boussinesq-Whitham (\ref{eqn:regBoussWhitham}) we note
\begin{align}
	\|\omega^2q^e  - \mathcal P\partial_x \left(  \: (\eta + \zeta)  (\mathcal Pq_x^e)\:  \right) \|_2 &\leq 
	\|\omega^2 q^e\|_2 + \| \mathcal P\partial_x \left(  \: (\eta + \zeta)  (\mathcal Pq_x^e)\:  \right) \|_2\\
	&\leq \frac{1}{\mu} \|k\tanh(\mu k)\hat q^e\|_2 +  \frac{1}{\mu}\|\eta+\zeta\|_\infty \|q_x^e\|_2,\\
	&\leq \frac{1}{\mu}\left(1 + \|\eta+\zeta\|_\infty \right)\|q_x^e\|_2,\label{eqn:bousswhitham:qx_estimate}
\end{align}
where we have assumed the data and true bottom boundary are bounded in the infinity norm (a reasonable supposition) and used the fact that $\mathcal P\partial_x$ corresponds to the multiplier $\tanh(\mu k)/\mu$ in Fourier space for regularised Boussinesq-Whitham. We also point out that although increasing the parameter $\mu$ can make the right-hand side of (\ref{eqn:error_zeta}) smaller, we have precisely the same scaling relationship with respect to $\mu$ on the left-hand side. Indeed setting $\mu\gg1$ makes the operator on the left-hand even more poorly conditioned. Informally speaking, we wish the right-hand side to be small, but the left-hand side to be `large', in the sense that the eigenvalues are bounded away from the origin. The scaling with respect to $\mu$ offers no advantage in this regard. Next we note the factor $(1 + \|\eta+\zeta\|_\infty )$ is always greater than unity. Hence there is a natural `amplification' of the error in $q_x$ in our proposed reconstruction method. Looking ahead, we will seek an estimate of $q_x^e$ that depends on an \textit{a priori} estimate of $\zeta^e$. Our error bound on $\zeta^e$ will only improve if the estimate of $q_x^e$ accounts for the amplification observed here.

The case for regularised Boussinesq (\ref{eqn:regBouss}) follows similarly. We have
\begin{align}
	\|\mathcal P\partial_x \big(\omega^2  q^e - \mathcal P\partial_x\left( (\mathcal P  q_x^e ) \: (\eta + \zeta) \right)  \big) \|_2 &\leq 
	\|\mathcal P\partial_x\omega^2 q^e\|_2  + \| \mathcal P^2\partial_x^2 \big(  (\eta + \zeta) \: (\mathcal P q_x^e) \big)\|_2\\
	&\leq \left\|\frac{k^2}{1+\mu^2k^2/2}\left(\frac{1+\mu^2k^2/6}{1+\mu^2k^2/2}\right)k \hat q^e\right\|_2 \nonumber \\
	&\quad 
	+ \left\| \frac{k^2}{(1+\mu^2k^2/2)^2} \mathcal{F}[(\eta+\zeta)(\mathcal P q_x^e)] \right\|_2,  \\
	&\leq\frac{1}{\mu^2}\left( 2  + \frac{1}{2}\|\eta+\zeta\|_\infty \right) \|q_x^e\|_2. \label{eqn:bouss:qx_estimate}
\end{align}
Apart from a different power of $\mu$ in the denominator, we have a similar amplification factor in front of the error in $q_x^e$. In fact the factors, for regularised Boussinesq and regularised Boussinesq-Whitham, are conservatively bounded above by $3$ assuming both $\|\eta\|_\infty$ and $\|\zeta\|_\infty$ are bounded above by $1$. This seems a reasonable assumption both from the perspective of having no islands but also from the point of view that these models approximate the full water-wave problem only in the small-amplitude shallow-water limit. For regularised Boussinesq, the scaling in $\mu$ on the left-hand side of (\ref{eqn:error_zeta}) is the same as the right-hand side. Arguably it is worse since on the left-hand we have a compact operator when $\mathcal P$ is given by (\ref{eqn:omegaPb}). Hence increasing $\mu$ only decreases the magnitude of the eigenvalues of the operator. Once again we conclude, it is the error in $q_x$ that must be made small.

The above discussion pertains to the case when $\eta_t,\eta,\tilde q_x$ were known at a single time. The estimates obtained, translate in a straightforward manner for the minimisation problem given in (\ref{eqn:bathymetry_min_alltime}). The only difference is we now demand estimates that are uniform in time. Hence we require $\tilde q_x$, the velocity as given by the observer model, to be sufficiently accurate over a period of time when $\eta$ and $\eta_t$ are also known. Here sufficiently accurate refers to the levels of accuracy required in the examples discussed in Section \ref{sec:bathy_reconstruct:examples}. We do not provide rigorous bounds on the required accuracy, which would entail an analysis of the spectrum of the operator on the left-hand side of (\ref{eqn:error_zeta}) as a \emph{function} of $\tilde q_x$.

{{\remark{It is worth mentioning the role of the operator $\mathcal P$ in obtaining the estimate (\ref{eqn:bouss:qx_estimate}). Note
\[\omega^2(k) = k^2 \left(\frac{1+\mu^2k^2/6}{1+\mu^2k^2/2}\right),\] represents an unbounded operator on $q$. Hence to interpret this expression as a bounded operator acting on $q_x$, the additional $\mathcal P\partial_x$ in the right-hand side of (\ref{eqn:error_zeta}) is precisely of the right form. In fact, had one considered a shallow-water model with a hyperbolic-tangent type regularisation on the nonlinear terms but retained the same $\omega^2$ as above, we would not be able to readily find a bound. Alternatively, had we considered a linear dispersion relation that included surface tension effects, we would require even greater smoothing from $\mathcal P$. This suggests the selection of the model is absolutely crucial for bottom-profile reconstruction. In particular it is the combination of the operators $\omega^2$ and $\mathcal P$ and their compatibility that is relevant.}}

\subsubsection{Error estimates for the observer problem}
We now focus our attention on the observer equations and the choice of observer parameters $\lambda,\nu$. As described in the reconstruction algorithm in Section \ref{sec:reconstruction_algo}, we initially guess the constant profile $\tilde\zeta(x) = \zeta_c$. This initial guess corresponds to an initial error given by $\zeta^e=\zeta_c-\zeta$. Consider the linear equation (\ref{eqn:observer_linear:zetaa}-\ref{eqn:observer_linear:zetab}) but now with an estimated bottom boundary:
\vspace{10pt}
\begin{align}\label{eqn:linear_observer:error_zetaa}
	 \eta_t^e &= \omega^2  q^e - \mathcal P\partial_x \big( (\mathcal P q_x^e) \: \zeta\big)-\lambda \eta^e - \mathcal P\partial_x \big( (\mathcal P q_x^e)\: \zeta^e \big) - \mathcal P\partial_x \big( (\mathcal P q_x)\: \zeta^e \big),\\
	 q_t^e &= -(1+\nu) \eta^e. \label{eqn:linear_observer:error_zetab}
\end{align}
Note the above equations are the error equations corresponding to the linearised version of (\ref{eqn:observer_NLa}-\ref{eqn:observer_NLb}) but with $\tilde \zeta$ in the place of $\zeta$. In terms of the observer problem, we seek to minimise the errors $q_x^e$ and $\eta^e$ in the presence of a forcing given by the last term of (\ref{eqn:linear_observer:error_zetaa}): $ \mathcal P\partial_x \big( (\mathcal P q_x)\: \zeta^e \big)$. 
As one might expect, due to the forcing term present when $\zeta$ is not known exactly, we cannot reduce the total error to zero. Moreover, the strength of this forcing is directly proportional to the error in $\zeta$. 

Equations (\ref{eqn:linear_observer:error_zetaa}-\ref{eqn:linear_observer:error_zetab}) can be combined to obtain a single equation for $q^e$
\begin{align}
	 q^e_{tt} + \lambda q^e_t + (1+\nu)\left(\omega^2q^e - \mathcal P\partial_x\big( (\mathcal Pq_x^e)\:  (\zeta + \zeta^e)\big) \right) = (1+\nu)\mathcal P\partial_x \left( (\mathcal Pq_x)\: \zeta^e\right).\label{eqn:q_err:linear}
\end{align}
Recall $\zeta^e = \zeta_c-\zeta$ where $\zeta_c$ is a constant. Hence the above is a constant-coefficient evolution equation for $q^e$ with a non-homogeneous term. This allows us to write an explicit solution using Fourier series. Indeed solutions to the homogeneous part are given in terms of exponential solutions of the kind $e^{ikx+pt}$ where $p$ satisfies
\[p^2 + \lambda p + (1+\nu)(\omega(k)^2 + \zeta_c\: k^2\mathcal P(k)^2) = 0.\]
Under the condition \[(1+\nu)> \left.\frac{\lambda^2}{4(\omega(k)^2 + \zeta_c\: k^2\mathcal P(k)^2)}\right|_{k=1},\] the roots are given by $p^{\pm} = -\lambda/2 \pm i \Omega_k$ where $\Omega_k$ is a non-zero real number for all $|k|>0$. See also the discussion in Section \ref{sec:observer:practice}. The general solution to (\ref{eqn:q_err:linear}) is then given by
\begin{align*}
q^e &= e^{-\frac{\lambda}{2}t}\sideset{}{'}\sum_{k=-\infty}^\infty e^{ikx}\left( e^{i\Omega_kt}\alpha_k  + e^{-i\Omega_kt}\beta_k\right) +e^{-\lambda t}\alpha_0 + \beta_0 \nonumber\\
&\quad + \int_0^t(1+\nu) e^{-\frac{\lambda}{2}(t-s)}\left(\sideset{}{'}\sum_{k=-\infty}^\infty e^{ikx}\frac{\sin(\Omega_k(t-s))}{\Omega_k}\mathcal{F}\left[\mathcal{P}\partial_x\big((\mathcal P q_x)\:\zeta^e\big)\right]_k\right)\:ds,
\end{align*}
where $\alpha_k,\beta_k$ depend on the Fourier coefficients of the initial data and the prime indicates the summation does not include $k=0$. We note that the error in $q$ is driven by a forcing proportional to $\zeta^e$, the error in our estimate of $\zeta$. We can now estimate $q_x^e$
\[\|q_x^e\|_2 \leq e^{-\frac{\lambda}{2}t}(\|\alpha_x\|_2 + \|\beta_x\|_2) + \left(\sup_t \|\mathcal{P}\partial_x^2\big((\mathcal P q_x)\:\zeta^e\big)\|_2\right) \int_0^t(1+\nu) e^{-\frac{\lambda}{2}(t-s)}ds, \]
and invoke Theorem \ref{thm:gronwalls}. Then for regularised Boussinesq, we expect as $t\to\infty$
\begin{align}
 \|q_x^e\|_2 \sim 2\frac{1+\nu}{\lambda}\sup_t \|\mathcal{P}\partial_x^2\big((\mathcal P q_x)\:\zeta^e\big)\|_2 \leq  \frac{4(1+\nu)}{\lambda\mu^2} \|\zeta^e\|_2 \:\sup_t \|\mathcal Pq_x\|_\infty.\label{eqn:observer_error:scale}
 \end{align}
 If $q_x\in L^2$ for all time, then as $\mathcal P$ is a smoothing operator, $\displaystyle\sup_t \|\mathcal P q_x\|_\infty$ is bounded. 

 The estimate (\ref{eqn:observer_error:scale}) implies we cannot choose $\lambda$ arbitrarily large to reduce the long-time error in $q_x$. Recall that $(1+\nu)$ is proportional to $\lambda^2$ in order to ensure the desired decay-rate. Hence increasing $\lambda$ only forces $\nu$ to be larger. We are then forced to make the combination $(1+\nu)/\lambda$ as small as possible so that the $\|q_x^e\|$ can be made sufficiently small.

We can make a similar argument to the one above for regularised Boussinesq-Whitham. However there is one technical obstacle. In this case, since $\mathcal P\partial_x$ is only bounded and not smoothing, we need to impose additional regularity for $\zeta^e$. For regularised Boussinesq-Whitham, the analogous estimate is given by
\begin{align}
 \|q_x^e\|_2 \sim 2\frac{1+\nu}{\mu\lambda}\sup_t \|\partial_x\big((\mathcal P q_x)\:\zeta^e\big)\|_2.
 \end{align}
Note we have a similar scaling with respect to $\nu$ and $\lambda$. Of course the above estimates are all based on an understanding of the linear equations  (\ref{eqn:linear_observer:error_zetaa}-\ref{eqn:linear_observer:error_zetab}). However, since the nonlinear terms in (\ref{eqn:observer_NLa}-\ref{eqn:observer_NLb}) are Lipschitz functions (for $\eta,q_x$ in $L^2$), we expect a similar estimate to hold for the nonlinear equations, if $\lambda$ is sufficiently large to overcome any possible growth in the solution. 

\subsubsection{Combining the estimates}
The key to reconstructing the bottom-boundary is to improve the error estimate. In other words, the final error in the bottom-profile $\zeta^e_{final}=\zeta^*-\zeta$, at the end of the algorithm described in Section \ref{sec:reconstruction_algo}, should be smaller than the initial error $\zeta^e_{init}=\zeta_c-\zeta$, where $\zeta_c$ is the initial guess for the bottom-profile. Consider the case of regularised Boussinesq. Upon combining the estimate (\ref{eqn:observer_error:scale}) with (\ref{eqn:bouss:qx_estimate}), we note the right-hand side of (\ref{eqn:error_zeta}) can be estimated as follows
\begin{align}
\| (\mathcal P\tilde q_x) \: \mathcal P^2\partial_x^2 \big( (\mathcal P\tilde q_x) \: \zeta^e_{final}\big)\|_2 &=
 \| (\mathcal P\tilde q_x)\: \mathcal P\partial_x \big(\omega^2  q^e - \mathcal P\partial_x\left( (\mathcal P  q_x^e ) \: (\eta + \zeta) \right)  \big)\|_2 \\
 &\leq	 \frac{12}{\mu^4} \frac{(1+\nu)}{\lambda} \|\zeta^e_{init}\|_2\: \|\mathcal P \tilde q_x\|_\infty \:\sup_t \|\mathcal Pq_x\|_\infty.
\end{align}
The error in the final estimated bottom profile will be small, if coefficient of $\|\zeta^e_{init}\|_2$ can be made as small as possible. For our simulations $|\tilde q_x|$ and $|q_x|$ are typically less than $1$. This is an outcome of our non-dimensional scaling and the initial conditions we used. Thus to ensure accurate reconstruction, we effectively need the combination $(1+\nu)/\lambda$ to be as small as possible. Although this conclusion only holds for regularised Boussinesq, we conjecture a similar estimate is true for regularised Boussinesq-Whitham. Our simulations indicate this is indeed the case.

\begin{figure}
\centering
\subcaptionbox{Error decay $\lambda = 1/100,\nu = -1 + \lambda^2$}[0.475\linewidth]{\includegraphics[width=0.475\textwidth]{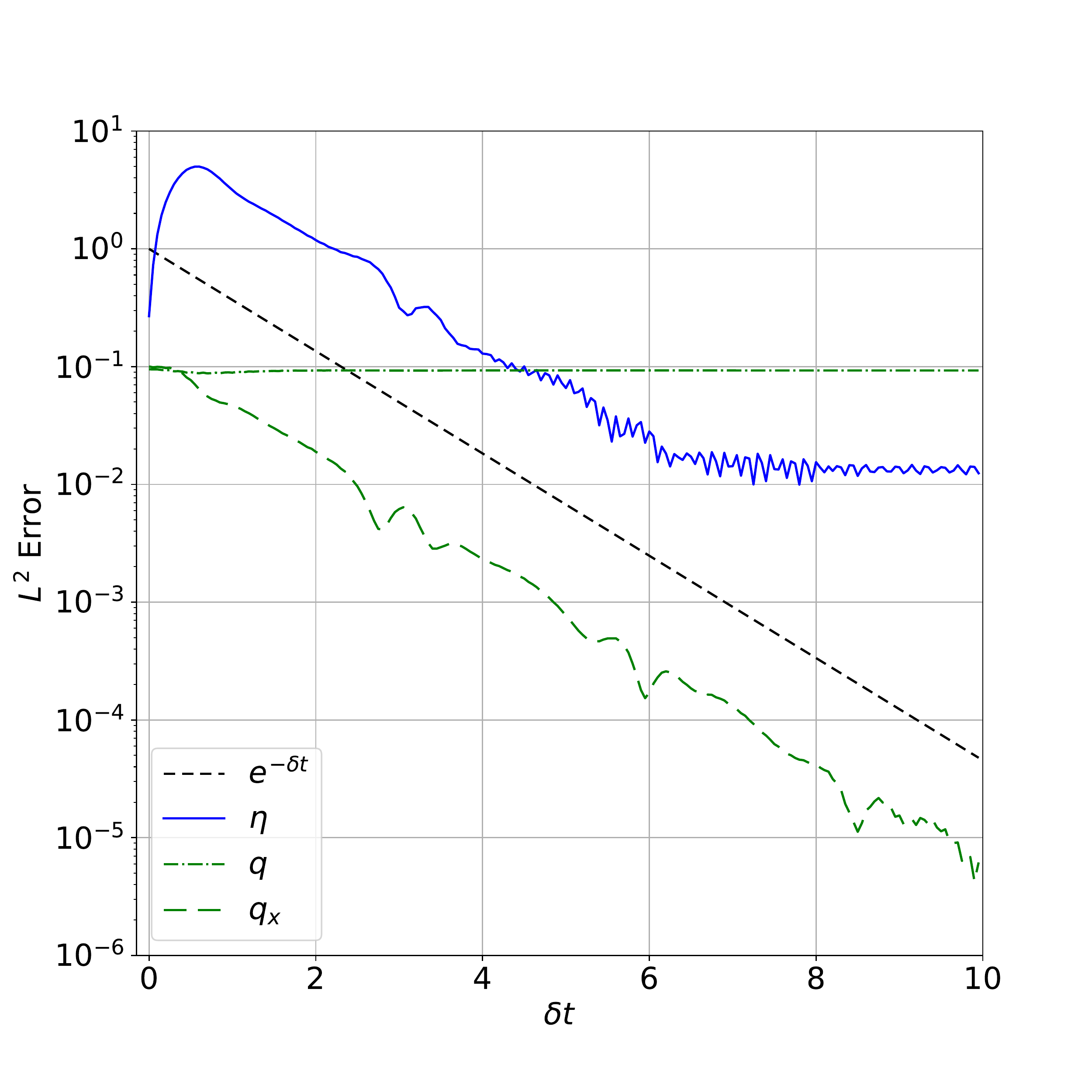}}
\subcaptionbox{Reconstruction }[0.475\linewidth]{\includegraphics[width=0.475\textwidth]{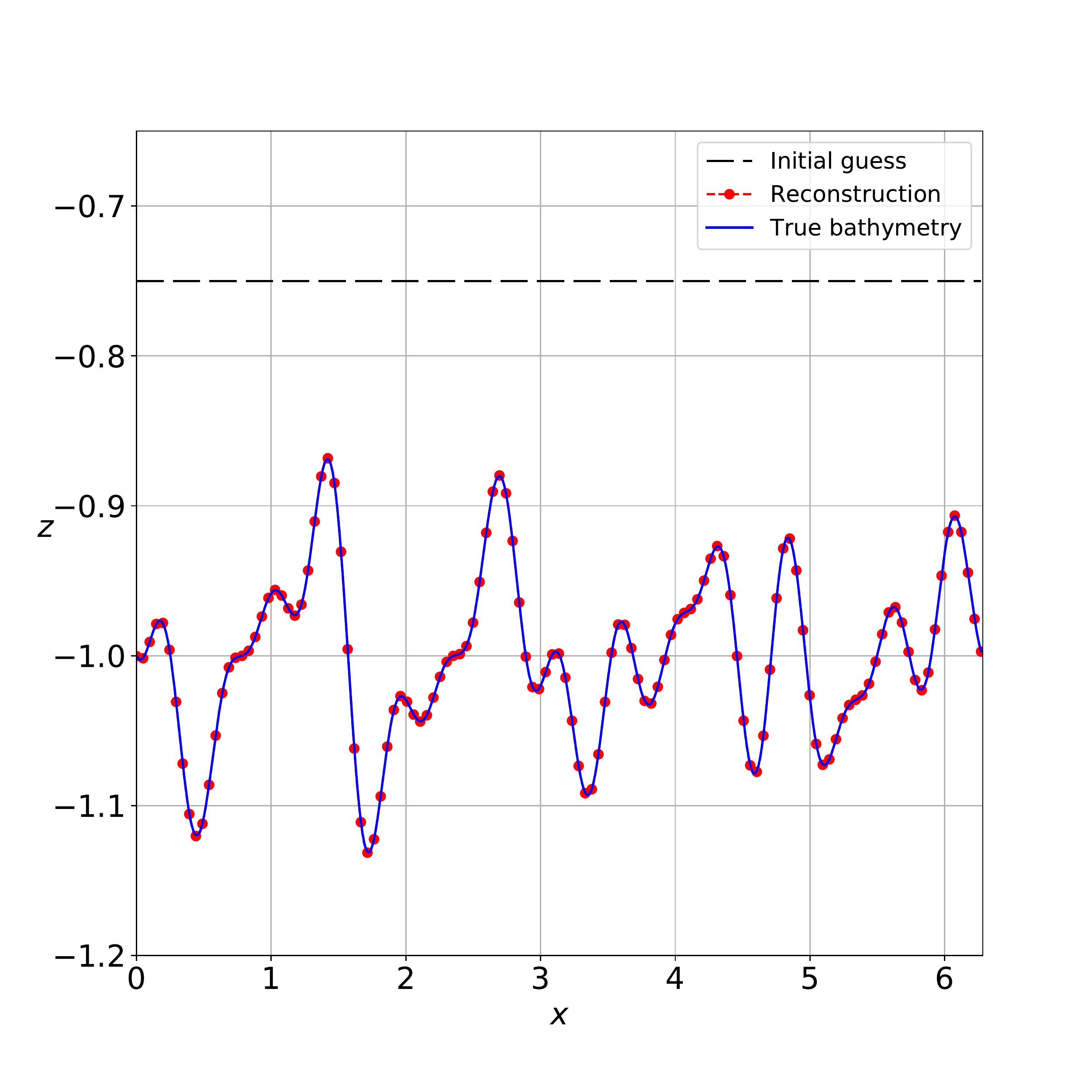}}
\caption{Reconstruction of \ref{fig:profile_wavy} using Regularised Boussinesq}\label{fig:simultaneous:Bouss:profile1}
\end{figure}

\begin{figure}
\centering
\subcaptionbox{Error decay $\lambda = 1/100,\nu = -1 + \lambda^2$}[0.475\linewidth]{\includegraphics[width=0.475\textwidth]{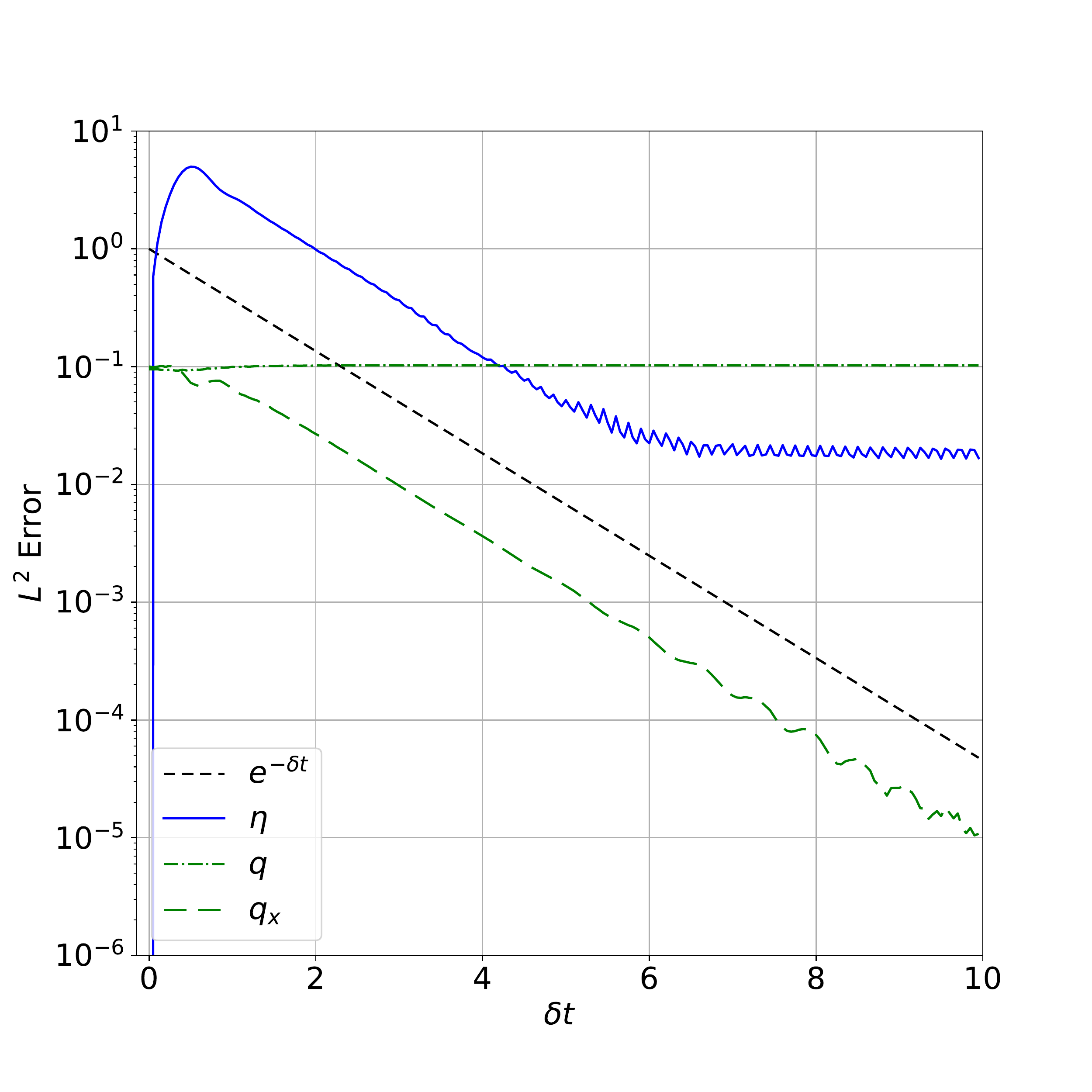}}
\subcaptionbox{Reconstruction }[0.475\linewidth]{\includegraphics[width=0.475\textwidth]{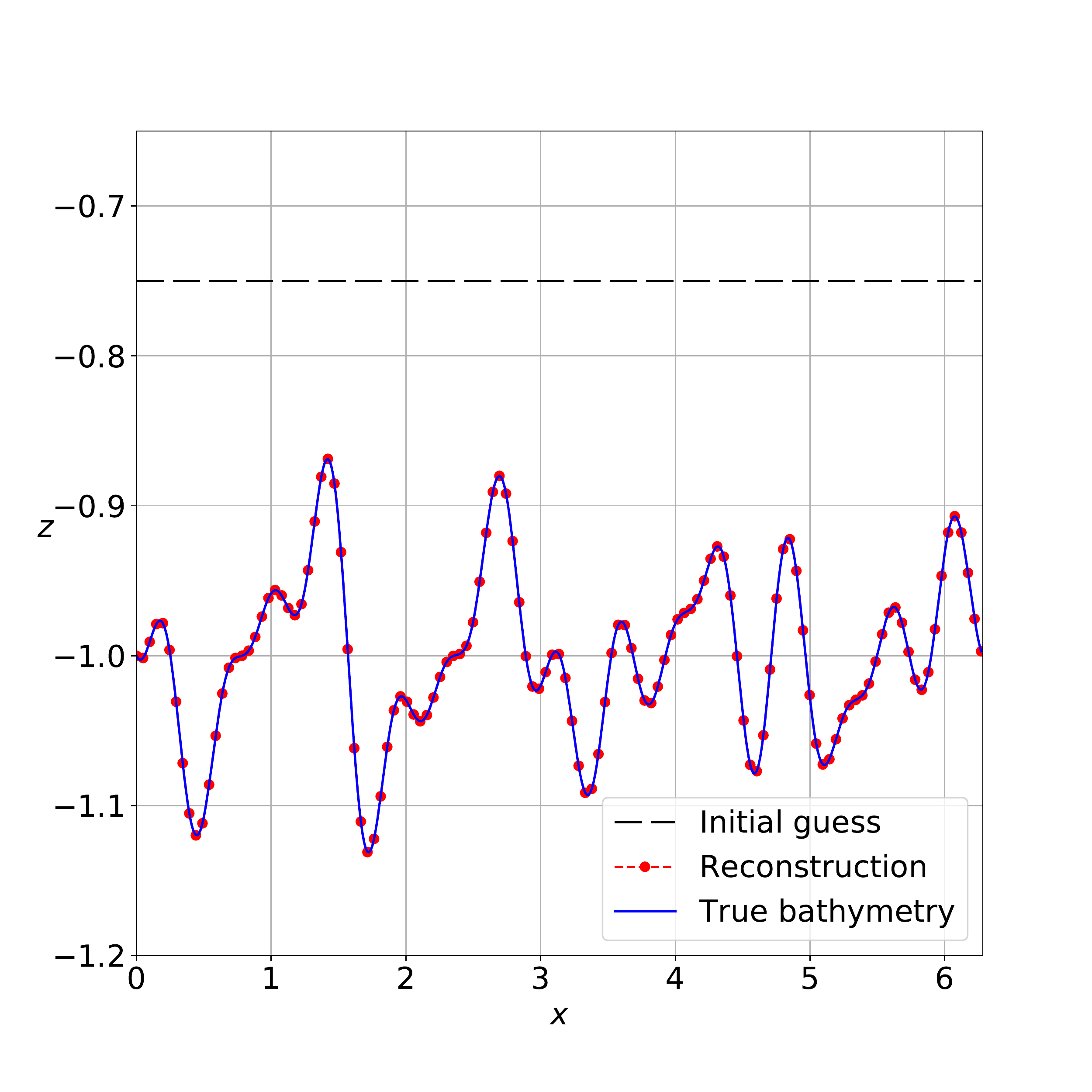}}
\caption{Reconstruction of \ref{fig:profile_wavy} using Regularised Boussinesq-Whitham}\label{fig:simultaneous:BoussWhitham:profile1}
\end{figure}

\begin{figure}
\centering
\subcaptionbox{Error decay $\lambda = 1/100,\nu = -1 + \lambda^2$}[0.475\linewidth]{\includegraphics[width=0.475\textwidth]{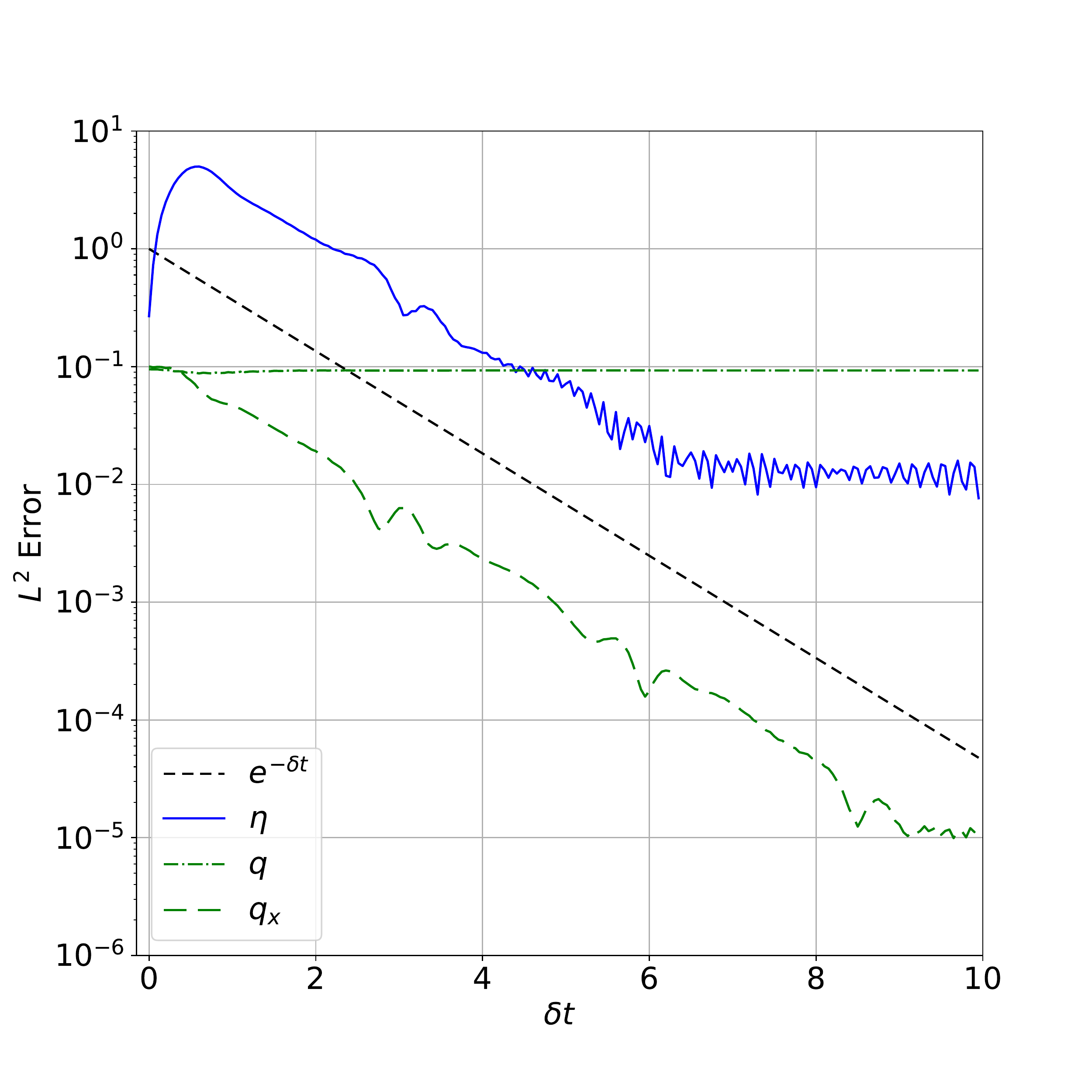}}
\subcaptionbox{Reconstruction }[0.475\linewidth]{\includegraphics[width=0.475\textwidth]{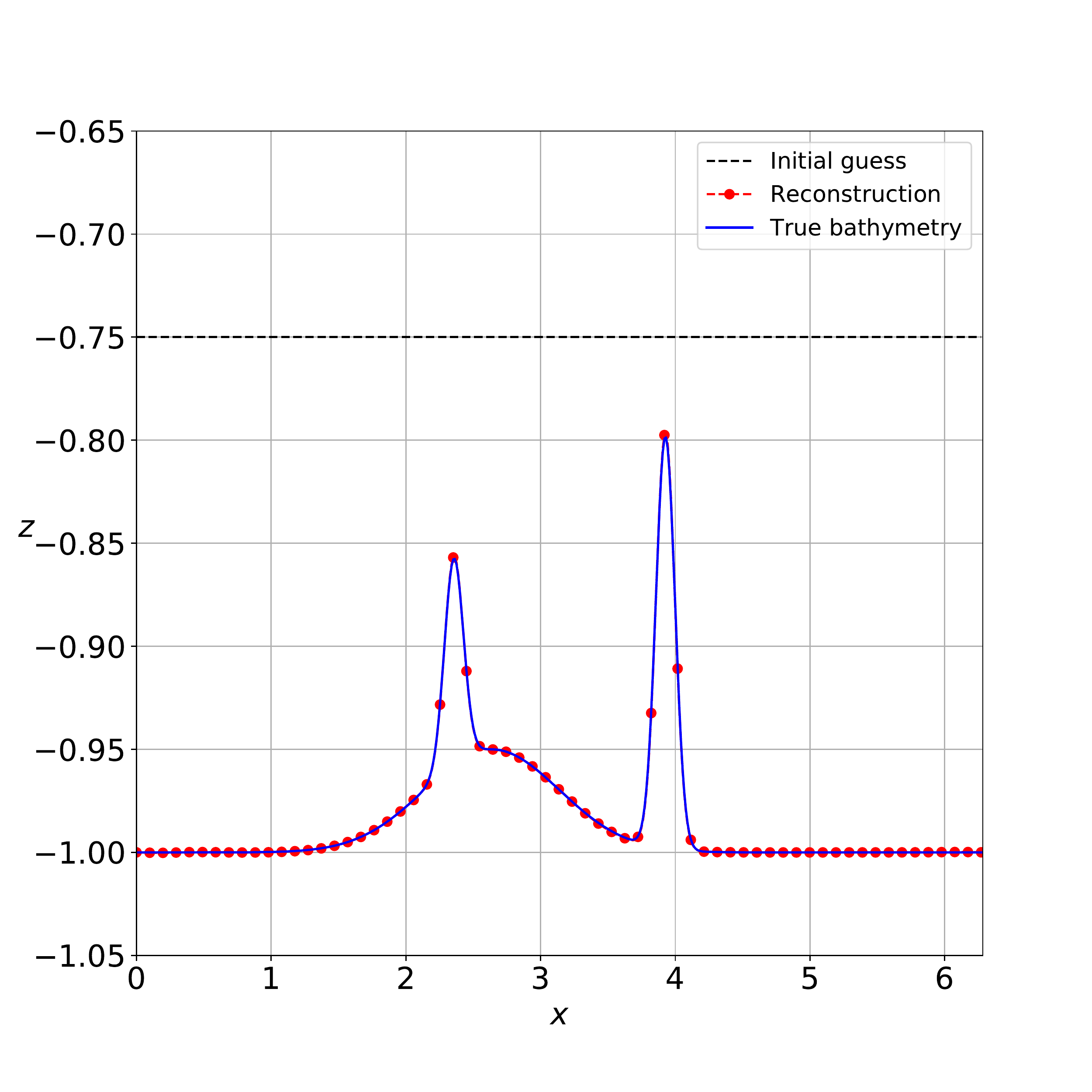}}
\caption{Reconstruction of \ref{fig:profile_isolated} using Regularised Boussinesq} \label{fig:simultaneous:Bouss:profile2}
\end{figure}

\begin{figure}
\centering
\subcaptionbox{Error decay $\lambda = 1/100,\nu = -1 + \lambda^2$}[0.475\linewidth]{\includegraphics[width=0.475\textwidth]{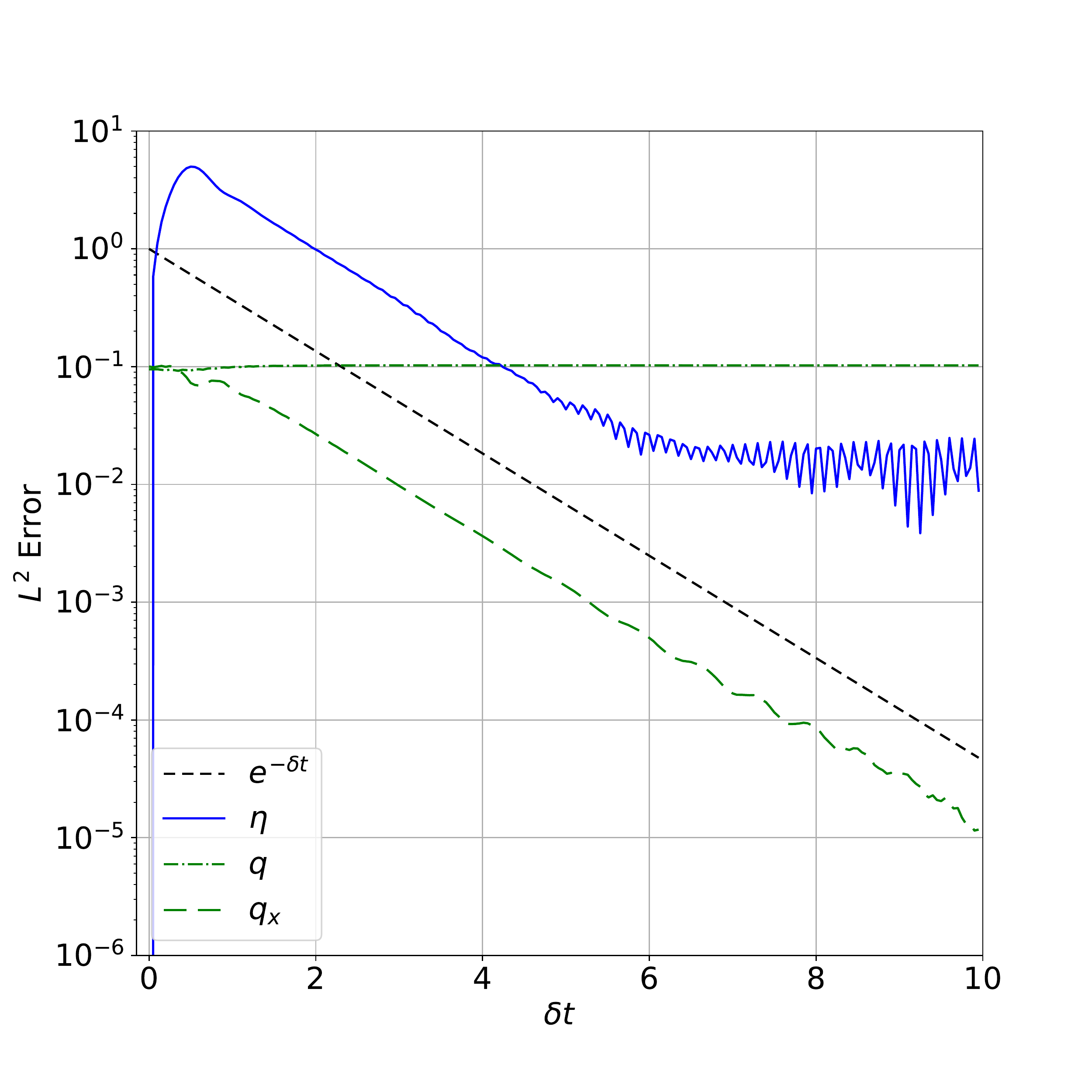}}
\subcaptionbox{Reconstruction }[0.475\linewidth]{\includegraphics[width=0.475\textwidth]{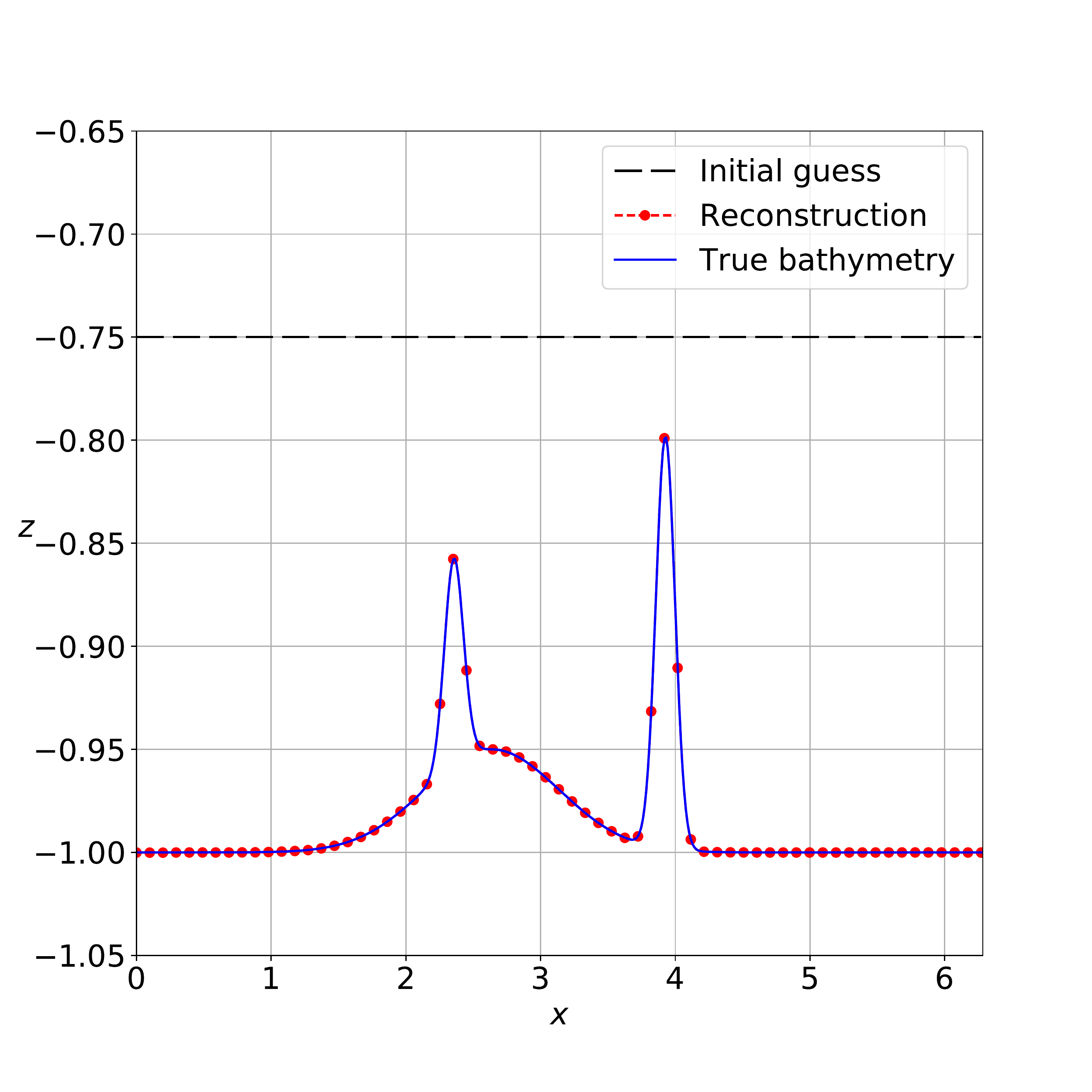}}
\caption{Reconstruction of \ref{fig:profile_isolated} using Regularised Boussinesq-Whitham}\label{fig:simultaneous:BoussWhitham:profile2}
\end{figure}

\subsection{Example reconstructions}

In Figures \ref{fig:simultaneous:Bouss:profile1} to \ref{fig:simultaneous:BoussWhitham:profile2} we show the result of applying the algorithm detailed above in the recovery of the two bottom profiles, Profile 1 and Profile 2 (see Figures \ref{fig:profile_wavy} and \ref{fig:profile_isolated}). The data necessary for reconstructing the bottom-boundary was obtained from a simulation of the model shallow-water equations (\ref{eqn:omegaPa}-\ref{eqn:omegaPb}) for either regularised Boussinesq or regularised Boussinesq-Whitham. In both cases we set the shallowness parameter $\mu=1$ which corresponds to an aspect ratio for the fluid $h/L \sim 0.16$. The initial condition for the shallow-water model is given by (\ref{init:conditions}) with $A = 0.0525$. The solution $\eta(x,t)$ of the shallow-water wave models was used to drive the observer model (\ref{eqn:observer_NLa}-\ref{eqn:observer_NLb}). The parameters for the numerical simulation of the observer were the same as those discussed in Section \ref{sec:experiments:observer}. The initial condition for the surface deviation in the observer model is the same as the one in (\ref{init:conditions}) but we set the initial velocity potential to zero. Our initial guess for the bottom boundary was $z=-0.75$ which corresponds to $\zeta_c = -0.25$. The observer parameters were $\lambda = 10^{-2}$ and $\nu = -1 + \lambda^2$ which ensured the linear decay rate of $\delta = \lambda/2$. 
The relative error $\mbox{E}_b$ (\ref{errors}) in the initial guess for the bottom-boundary was approximately $25\%$ or $23\%$ (for Profile 1 and Profile 2 respectively) which corresponds to an error of over $430\%$ in terms of $\mbox{E}_p$ (\ref{errors}). For either profile, the values for $\mbox{E}_b$ for recovered bottom-boundary was $2\times 10^{-4}$ for regularised Boussinesq and $6\times 10^{-5}$ for regularised Boussinesq-Whitham. In terms of $\mbox{E}_p$, these are $99.5\%$ and $99.9\%$ accurate reconstructions. We did not observe the error in the initial guess of the bottom-profile to dramatically impact the recovery, so long as the error in the estimated velocity $\tilde q_x$ could be reduced to a sufficiently low value. A relative error of $10^{-4}$ in the estimated velocity (measured in the $L^2-$ norm) was sufficient in all our examples. As evident in Figures \ref{fig:simultaneous:Bouss:profile1} to \ref{fig:simultaneous:BoussWhitham:profile2}, the error in the velocity follows the predicted linear rate whereas the error in the surface deviation saturates around $10^{-2}$. The surface deviation from the observer model $\tilde \eta$, is not needed for the reconstruction algorithm and hence the higher level of error is not a cause of concern. The observer model was run till a nondimensional time of $2000$ units. Note the horizontal axis in the error plots of Figures \ref{fig:simultaneous:Bouss:profile1} to \ref{fig:simultaneous:BoussWhitham:profile2} is given in terms of $\delta t$ where $\delta = 0.005$. Our proposed algorithm successfully recovers the bottom-boundary with both shallow-water wave models and for different bottom profiles starting from a relatively inaccurate initial guess.

\section{Summary and future work}\label{sec:summary}
We have shown that the simultaneous recovery of both velocities and bottom-boundary from only measurements of the surface deviation, in the context of dispersive shallow-water models, is possible. The motivation to consider bottom-boundary detection in shallow water comes from prior work \cite{vasan2013inverse}. The aspect ratios $h/L$ for which we are able to recover the bottom profiles are larger,  by an order of magnitude, than what was previously considered in \cite{vasan2013inverse}. This suggests there is some advantage to our prescription. We have also shown that velocimetry, the estimation of velocities, in Boussinesq-type shallow-water models is entirely feasible when given an accurate bottom-boundary profile. Moreover, using the observer framework, velocities can be accurately estimated even when the bottom-boundary profile is \textit{not} known. Unlike traditional observers, for velocimetry with inaccurate bottom-boundary profiles, we cannot use the convergence of the measured state-variable (here $\eta$) to infer convergence of the full state-vector. Instead we employed the expected decay rate to ascertain when the velocity was sufficiently accurate. Lastly, we have shown that the ability to recover the bottom-boundary and velocity is not restricted to a particular shallow-water model. Indeed it seems recovery is possible for a class of models depending on the pseudo-differential operators $\omega^2,\mathcal P$.

Despite our success, we emphasise the ocean-depth measurement is a delicate process that involves carefully setting the parameters for the observer problem. The requirements for the simultaneous recovery of velocities and bottom profiles are somewhat contradictory. The nonlinear observer problem requires $\lambda$ to be sufficiently large to dominate any possible growth in the error due to the nonlinearity. On the other hand, recovering bathymetric features requires accurate $\tilde q_x$ which demands $\lambda$ be taken as small as possible. Our numerical experiments indicate there is a parameter regime for which one may accurately reconstruct the bottom-boundary profile. In a future paper we will explore the possibility for time-dependent observer parameters $\lambda,\nu$ that slowly relax from larger to smaller values and whether this provides any advantage. 

All our simulations were conducted for one spatial variable, however the extension of the shallow-water model to two horizontal dimensions is straightforward. Indeed 
\begin{align}
    \eta_t &= \omega^2(-i\nabla)q - \mathcal{P}(-i\nabla) \nabla\cdot \left[ (\eta+ \zeta)\mathcal{P}(-i \nabla)\nabla q\right],\\ 
	q_t &= -\eta - \frac 1 2 \left(\mathcal{P}(-i\nabla)\nabla q \right)^2,
\end{align}
are the $2D$ versions, where $\omega^2$ and $\mathcal P$ are defined as before but as functions of $|k|$. The observer equations may be similarly rewritten in two spatial dimensions. Most of our arguments can be adapted to the $2D$ case without difficulty. We leave the full investigation of the $2D$ ocean-depth measurement problem for a future work.

The simultaneous recovery of velocities and bottom-boundaries necessitated very long time simulations for the observer model. One might suspect this requires an equally long data record for $\eta(x,t)$. However this is not necessarily the case. The shallow-water models are Hamiltonian and thus time-reversible. This permits us to assimilate the data into the observer model in both the forward and backward directions of time, taking care to re-index the data record. This is the principle underlying back and forth nudging \cite{AurouxBlum08,auroux2005back}. To be precise, suppose we only had measurements of $\eta$ for a finite duration of time $[0,T]$. We can run the observer model forward in time up to time $T$ and then run the model backward in time using the final state of the observer as the new initial condition. In the backward-run, the data $\eta(x,t)$ is reversed in time. 

For real-world applications, it is important to adapt our reconstruction algorithm to the case with non-periodic data $\eta(x,t)$. This would represent a significant improvement upon the problem as described in the current manuscript. When the domain is not periodic, $\int\eta$ is not necessarily a constant and thus the system involves mass flux across the lateral boundaries. The successful modelling of this scenario would entail a re-derivation of the associated DNO and the Hamiltonian formulation, and subsequently obtaining the relevant approximations.

The current work did not consider the addition of noise in the measurements. As seen in other observer problems \cite{auroux-bonnabel11}, adding noise to the measurements in the observer model will likely saturate the errors in $\tilde q_x$ to non-zero values at the level of the noise, though the observer formulation itself is easy to adapt. Since a low relative-error in the estimated velocity was crucial to obtaining accurate bottom-profile reconstructions, noise will have a significant impact on the accuracy of the reconstructed profiles. As in \cite{vasan2013inverse}, computing the time-derivative of $\eta$ using finite differences is straightforward and not the main source of error. Instead the accuracy of the estimated velocities determines the ultimate success of the reconstruction algorithm. 

And finally we conclude by emphasising the role implicitly played by model selection in the entire reconstruction process. In general, the regularised Boussinesq model has a number of favourable aspects that enabled us to estimate quantities of interest. However, the drawback of this model was the rapid of decay of eigenvalues for the linear operator in (\ref{eqn:bathymetry_EL:all_time}). The regularised Boussinesq-Whitham model did not suffer from such a rapid decay however, we were not able to verify some of the estimates or required additional regularity assumptions on the bottom-profile. The reduction to the shallow-water regime was also important in allowing us to design the observer and the reconstruction algorithm in a relatively simple manner. Ultimately, the choice of model is dictated by our understanding and interpretation of the data, as well as our assumptions on the true bottom-boundary. However, for the shallow-water model to actually model the full water-wave equations, these two factors may not be independent of each other \cite{lannes2013water}. Our work shows that the mathematical properties of these equations, specifically the interdependence of the dispersion relation $\omega^2$ and regularisation $\mathcal P$, implies some models may be preferable to others. We hope the combination of simple estimates and numerical simulations reported in this manuscript, afforded the reader some insight into the nature of this difficult inverse problem.




\section*{Acknowledgements}
The authors thank Amit Apte and Mythily Ramaswamy for all their helpful comments and suggestions. The authors acknowledge support of the Department of Atomic Energy, Government of India, under project no. RTI$4001$. VV acknowledges support through the SERB MATRICS Grant (MTR/2019/000609) from the Science and Engineering Research Board (SERB), Department of Science and Technology, Government of India. This work has been supported by the French government, through the $\mbox{UCAJEDI}$ Investments in the Future project (UCA-ICTS joint program) managed by the National Research Agency (ANR) with the reference number ANR-15-IDEX-01. The authors thank the Indo-French Centre for Applied Mathematics (IFCAM) for financial support under the ``Observers for coupled models and parameter estimation'' project.

\appendix{}

\titleformat{\section}[runin]
{\color{RoyalBlue}\bfseries}{\color{RoyalBlue}\appendixname~\thesection .}{1ex}{}[.]

\section{Minimisation problems in the shallow water regime}\label{app:VD_shallow}
Fontelos et al. \cite{fontelos2017bottom} minimised the functional \[F(\zeta) = \int_0^{2\pi} (\eta_t - G(\eta,\zeta)q )^2dx, \] to obtain the bottom-profile $\zeta$, given $\eta,\eta_t,q$ at one time instant. They showed the minimiser of this functional corresponded to the true bottom-profile. In \cite{vasan2013inverse}, the authors employed a slightly different methodology to obtain a similar minimisation problem. They too defined a functional dependent on the bottom boundary, but in this case they derived an expression for the Neumann condition at some bottom surface $\zeta$ (dependent on surface quantities $\eta,\eta_t,q_x$) and sought to minimise this quantity, thereby enforcing the no-normal flow condition. We now show through formal calculations that, in the shallow-water regime, both approaches lead to the same minimisation problem (\ref{eqn:bathymetry_min}).

To derive asymptotic expressions for the Neumann condition at the bottom boundary, that are consistent with the models introduced in the previous section, we once again employ the global relation introduced by \cite{afm}. However our focus will shift slightly. Specifically, we assume \emph{both} conditions at the surface $z=\eta$ are known and \emph{both} boundary conditions at $z=-h-\zeta$ are unknown. The global relations are given by
\begin{align}
	\int e^{ikx}\left\{ in^t \mathrm{C}_k(\eta+1) + q_x\mathrm{S}_k(\eta+1) + Q_x\mathrm{S}_k(\zeta) + in^b \mathrm{C}_k(\zeta) \right\}dx &= 0,\\
	\int e^{ikx}\left\{ in^t \mathrm{S}_k(\eta+1) + q_x\mathrm{C}_k(\eta+1) - Q_x\mathrm{C}_k(\zeta) - in^b \mathrm{S}_k(\zeta) \right\}dx &= 0,
\end{align}
where $\mathrm{C}_k(f) = \cosh(\mu k f)$ and $\mathrm{S}_k(f)=\sinh(\mu k f)$ and $n^t,n^b$ denote the Neumann condition at the top and bottom respectively. $Q_x,q_x$ denote the same quantities as before. An expression for the bottom Neumann condition $n^b$ consistent with the choice (\ref{eqn:regBouss}) is obtained via the substitution 
\[\mathrm{C}_k(\eta + 1) \to  1 + \frac{(\mu k)^2}{2},\quad \mathrm{S}_k(\eta + 1) \to \mu k + \mu^3 k\eta + \frac{(\mu k)^3}{6}\]
\[ \mathrm{C}_k(\zeta)\to 1,\quad \mathrm{S}_k(\zeta)\to \mu^3 k\zeta. \]
This leads to $Q_x = q_x + \mbox{h.o.t.}$ and 
\begin{align}
	-n^b = \left(1-\frac{\mu^2 \partial_x^2}{2}\right)n^t + q_{xx} - \frac{\mu^2}{6}q_{xxxx} + \partial_x\left( (\eta+\zeta)q_x\right) + \mbox{h.o.t.}
\end{align}
Evidently, requiring the Neumann condition $n^b$ to vanish to the same order as regularised Boussinesq, for some $\zeta$, is equivalent to requiring 
\begin{align}
	{\int (n^t - \omega^2 q + \mathcal P \partial_x((\eta+\zeta)\mathcal P q_x) )^2 \: dx=0,
	}
\end{align}
with the choice (\ref{eqn:regBouss}). In exactly the same manner, the equation for the bottom Neumann-condition consistent with the approximations that lead to regularised Boussinesq-Whitham, is given by the same expressions as above but with the choice (\ref{eqn:regBoussWhitham}). Thus in the shallow-water regime, the approaches of \cite{fontelos2017bottom} and \cite{vasan2013inverse} are formally equivalent.

\section{Energy for the linear observer problem}\label{app:energy}
Consider an equivalent form of (\ref{eqn:observer_linear:zetaa}-\ref{eqn:observer_linear:zetab})
given by
\begin{align}\label{eqn:eta_linear:error}
	\eta_{tt}^e + \lambda \eta_t^e + (1+\nu)(\omega^2 \eta^e -\mathcal P \partial_x \big( \zeta \: (\mathcal P \eta_x^e) \big) = 0,
\end{align}
from which we deduce the energy relation
\begin{align}
	\frac 1 2 \partial_t\left( \int (\eta_t^e)^2 +(1+\nu) \int\left( \eta^e\: \omega^2 \eta^e + \zeta (\mathcal P\eta_x^e)^2 \right)\right)\: dx = -\lambda \int (\eta_t^e)^2 \:dx.
\end{align}
One might reasonably suppose that the bottom boundary does not penetrate the free surface, at least for the trivial solution. This is sometimes known as the non-cavitation or no-island condition. For the full nonlinear problem one demands the fluid height $1+\eta+\zeta>0$. The analogue for the linear equation here is $1+\zeta>0$. Then we have a positive energy if
\begin{align}
	\int (\eta^e\: \omega^2 \eta^e - (\mathcal P\eta_x^e)^2) \:dx &\geq0.
\end{align}
At this stage we must treat our two models independently. We first consider regularised Boussinesq (\ref{eqn:regBouss}) which leads to 
\begin{align}
	\int (\eta^e\: \omega^2 \eta^e - (\mathcal P\eta_x^e)^2) \:dx &= \sum_k \left( k^2\left(\frac{1+\mu^2k^2/6}{1+\mu^2k^2/2}\right) |\hat\eta^e_k|^2 - \frac{k^2|\hat\eta_k^e|^2}{(1+\mu^2k^2/2)^2}  \right),\\
	&= \sum_{k}\left(k^2|\hat\eta_k^e|^2\left(  \frac{1+\mu^2k^2/6}{1+\mu^2k^2/2} - \frac{1}{(1+\mu^2k^2/2)^2} \right) \right)\geq 0.
\end{align}
Likewise for regularised Boussinesq-Whitham (\ref{eqn:regBoussWhitham}) we have
\begin{align}
	\int (\eta^e\: \omega^2 \eta^e - (\mathcal P\eta_x^e)^2) \:dx &= \sum |\hat \eta_k^e|^2 \left(k\frac{\tanh(\mu k)}{\mu} - \frac{\tanh^2(\mu k)}{\mu^2}  \right),\\
	&= \sum_k k\frac{\tanh(\mu k)}{\mu}|\hat \eta_k^e|^2\left( 1 - \frac{\tanh(\mu k)}{\mu k} \right) \geq 0.
\end{align}
Hence the energy in both cases is positive so long as $\eta^e$ is not a constant in space. If we further assume the mean of $\eta^e$ remains zero for all time, \emph{i.e.} there is no error in the mean value of the surface deviation, then the potential energy may be bounded below. For the problem with periodic boundary conditions, the discrete nature of the spectrum is to our advantage. For regularised Boussinesq we have 
\begin{align}
	\int (\eta^e\: \omega^2 \eta^e - (\mathcal P\eta_x^e)^2) \:dx &\geq C_1(\mu)\int(\eta^e_x)^2,\quad {C_1(\mu) = \frac{\mu^2(2/3 + \mu^2/12)}{(1+\mu^2/2)^2} },
\end{align}
whereas for regularised Boussinesq-Whitham we have
\begin{align}
	\int (\eta^e\: \omega^2 \eta^e - (\mathcal P\eta_x^e)^2) \:dx &\geq \frac{C_2(\mu)}{\mu} \sum_{k\neq 0} k\tanh(\mu k)|\hat \eta_k^e|^2,\quad {C_2(\mu) = 1 - \frac{\tanh(\mu)}{\mu}},
\end{align}

\noindent Our purpose in considering the energy relation, is to investigate when the error will decrease monotonically, at least for the linear observer problem. Indeed the solutions which prevent monotonic decrease in the error are the steady solutions: $\eta_t^e = 0$ for all time. However (\ref{eqn:eta_linear:error}) then implies
\[\omega^2\eta^e - \mathcal P \partial_x \big(\zeta \:(\mathcal P \eta_x^e)\big) = 0 \quad \Rightarrow \int \eta^e\:\left(\omega^2\eta^e - \mathcal P \partial_x (\zeta \mathcal P \eta_x^e) \right)\: dx  = 0. \]
From the lower bound obtained above we conclude that only constant $\eta^e$ (in time and space) prevent any decay in the error.  But once again, the average of (\ref{eqn:observer_linear:zetaa}-\ref{eqn:observer_linear:zetab}) in the $x-$direction indicates it is the mean-mode of $q^e$ that does not decay to zero. The zero mode of $\eta^e$ decays exponentially.  This is true for both regularised Boussinesq and regularised Boussinesq-Whitham.

\bibliographystyle{amsplain}

\bibliography{refs.bib}

\end{document}